\documentclass[12pt,english,oneside,thmsb]{amsart}
\usepackage[T1]{fontenc}
\usepackage[latin1]{inputenc}
\usepackage{setspace}
\usepackage{amssymb}

\makeatletter

\newcommand{\noun}[1]{\textsc{#1}}

\numberwithin{equation}{section} 
\numberwithin{figure}{section} 
  \@ifundefined{theoremstyle}{\usepackage{amsthm}}{}
  \theoremstyle{plain}
  \newtheorem{thm}{Theorem}[section]
  \theoremstyle{remark}
  \newtheorem*{rem*}{Remark}
  \theoremstyle{plain}
  \newtheorem{lem}[thm]{Lemma}
  \theoremstyle{remark}
  \newtheorem{rem}[thm]{Remark}
  \theoremstyle{plain}
  \newtheorem{cor}[thm]{Corollary}

\usepackage{babel}
\makeatother
\def\1{\hbox{\rm 1\hskip -3pt I}}

\makeatother

\usepackage{babel}

\begin{document}

\author{Jean-Michel Loubes$^{\flat}$ and Anne-Françoise Yao $^{\S}$ }

\thanks{$\S$ Corresponding author : University Aix-Marseille 2, Campus de
Luminy, case 901, 13288 Marseille cedex 09. anne-francoise.yao@univmed.fr}

\thanks{$^{\flat}$ University Toulouse 3, Institut de Mathématiques de Toulouse.}

\title{Kernel Inverse Regression for spatial random fields.}

\begin{abstract}
In this paper, we propose a dimension reduction model for spatially
dependent variables. Namely, we investigate an extension of the \emph{inverse
regression} method under strong mixing condition. This method is based
on estimation of the matrix of covariance of the expectation of the
explanatory given the dependent variable, called the \emph{inverse
regression}. Then, we study, under strong mixing condition, the weak
and strong consistency of this estimate, using a kernel estimate of
the \emph{inverse regression}. We provide the asymptotic behaviour
of this estimate. A spatial predictor based on this dimension reduction
approach is also proposed. This latter appears as an alternative to
the spatial non-parametric predictor.
\end{abstract}
\maketitle
\noun{Keywords:} Kernel estimator; Spatial regression; Random fields;
Strong mixing coefficient; Dimension reduction; Inverse Regression.

\section{Introduction}

Spatial statistics includes any techniques which study phenomenons
observed on spatial subset $S$ of $\mathbb{R}^{N},\, N\geq2$ (generally,
$N=2$ or $N=3$). The set $S$ can be discret, continuous or the
set of realization of a point process. Such techniques have various
applications in several domains such as soil science, geology, oceanography,
econometrics, epidemiology, forestry and many others (see for example
\cite{riplet}, \cite{cressie} or \cite{guyon} for exposition, methods
and applications). 

Most often, spatial data are dependents and any spatial model must
be able to handle this aspect. The novelty of this dependency unlike
the time-dependency, is the lack of order relation. In fact, notions
of past, present and futur does not exist in space and this property
gives great flexibility in spatial modelling. 

In the case of \emph{spatial regression} that interests us, there
is an abundant literature on \emph{parametric models}. We refer for
example to the spatial regression models with correlated errors often
used in economics (see e.g. Anselin and Florax \cite{anselin}, Anselin
and Bera \cite{anselinbera98}, Song and Lee \cite{Song-Lee08}) or
to the spatial Generalized Linear Model (GLM) study in Diggle et al.
\cite{Diggleetal98} and Zhang \cite{Zang02}. Recall also the spatial
Poisson regression methods which have been proposed for epidemiological
data (see for example Diggle \cite{Diggle03} or Diggle et al \cite{Diggleetal98}).

Unlike the parametric case, the spatial regression on nonparametric
setting have been studied by a few paper: quote for example Biau and
Cadre \cite{biau-cadre}, Lu and Chen \cite{lu-chen04}, Hallin et
al. \cite{Hallin-lu-tran04b}, Carbon et al. \cite{carbonfrancqtran07},
Tran and Yakowitz \cite{tran-yako} and Dabo-Niang and Yao \cite{dabo-yao07}.
Their results show that, as in the \emph{i.i.d.} case, the spatial
nonparametric estimator of the regression function is penalized by
the dimension of the regressor. This is the spatial counterpart of
the well-known problem called {}``\emph{the curse of dimensionality}''.
Recall that dimension reduction methods are classically used to overcome
this issue. Observing an \emph{i.i.d.} sample $Z_{i}=(X_{i},Y_{i})$
the aim is to estimate the regression function $m(x)=\mathbf{E}(Y|X=x)$.
In the dimension reduction framework, one assumes that there exist
$\Phi$ an orthonormal matrix $d\times D$, with $D$ as small as
possible, and $g\,:\mathbb{R}^{D}\to\mathbb{R}$, an unknown function
such that the function $m(.)$ can be written as \begin{equation}
m(x)=g(\Phi\,.X).\label{eq:mod}\end{equation}
 Model~(\ref{eq:mod}) conveys the idea that {}``less information
on $X$'' , $\Phi\,.X$; provides as much information on $m(.)$
as $X$. The function $g$ is the regression function of $Y$ given
the $D$ dimensional vector $\Phi.X$. Estimating the matrix $\Phi$
and then the function $g$ (by nonparametric methods) provides an
estimator which converges faster than the initial nonparametric estimator.
The operator $\Phi$ is unique under orthogonal transformation. An
estimation of this latter is done through an estimation of his range
$\textrm{Im}(\Phi^{T})$ (where $\Phi^{T}$ is the transpose of $\Phi$)
called \emph{Effective Dimensional Reduction space} (EDR).\\
 \indent Various methods for dimension reduction exist in the literature
for \emph{i.i.d} observations. For example we refer to the multiple
linear regression, the generalized linear model (GLM) in \cite{brillinger},
the additive models (see e.g. Hastie and Tibshirani \cite{Hastie-Tibshi86})
deal with methods based on estimation of the gradient of the regression
function $m(.)$ developped in for example in \cite{Hirst-Judis-Polz-spok}
or \cite{xia-tong-li-zhu}. 

\vskip .1in 

In this paper, we focus on the \emph{inverse regression} method, proposed
by Li \cite{Li}: if \emph{$X$ is such that for all vector $b$ in
$\mathbb{R}^{d}$, there exists a vector $B$ of $\mathbb{R}^{D}$
such that $\mathbf{E}(b^{T}X|\Phi.X)=B^{T}(\Phi.X)$} (this latter
condition is satisfied as soon as $X$ is elliptically distributed),
then, if $\Sigma$ denotes the variance of $X$, the space $\mbox{Im}(\Sigma^{-1}\mathbf{var}(\mathbf{E}(X|Y))$
is included into the \emph{EDR space}. Moreover, the two spaces coincide
if the matrix $\Sigma^{-1}\mathbf{var}(\mathbf{E}(X|Y))$ is of full
rank. Hence, the estimation of the \emph{EDR space} is essentially
based on the estimation of the covariance matrix of the \emph{inverse
regression} \textbf{$\mathbf{E}(X|Y)$} and $\Sigma$ which is estimated
by using a classical empirical estimator. In his initial version,
Li suggested an estimator based on the regressogram estimate of $\mathbf{E}(X|Y)$
but drawbacks of the regressogram lead other authors to suggest alternatives
based on the nonparametric estimation of \textbf{$\mathbf{E}X|Y$},
see for instance \cite{Hsing99} or \cite{Fang-Zhu96} which enable
to recover the optimal rate of convergence in $\sqrt{n}$. 

This work is motivated by the fact that to our knowledge, there is
no \emph{inverse regression} method estimation for spatially dependent
data under strong mixing condition. Note however that a dimension
reduction method for supervised motion segmentation based on spatial-frequential
analysis called \emph{Dynamic Sliced Inverse Regression} (DSIR) has
been proposed by Wu and Lu \cite{Wu-Lu-04}. We propose here a spatial
counterpart of the estimating method of \cite{Fang-Zhu96} which uses
kernel estimation of \textbf{$\mathbf{E}X|Y$}. Other methods based
on other spatial estimators of \textbf{$\mathbf{E}X|Y$} will be the
subject of futher investigation. 

\vskip .1in 

As any spatial model, a spatial dimension reduction model must take
into account spatial dependency. In this work, we focus on\emph{ }an
estimation on model (\ref{eq:mod}) for spatial dependent data under
strong mixing conditions. The spatial kernel regression estimation
of \textbf{$\mathbf{E}X|Y$} being studied in \cite{biau-cadre,carbon-tran-wu,carbonfrancqtran07}. 

\vskip .1in 

An important problem in spatial modelling is that of spatial prediction.
The aim being reconstruction of a random field over some domain from
a set of observed values. It is such a problem that interest us in
the last part of this paper. More precisely, we will use the properties
of the \emph{inverse regression} method to build a \emph{dimension
reduction predictor} which corresponds to the \emph{nonparametric
predictor} of \cite{biau-cadre}. It is an interesting alternative
to parametric predictor methods such as the \emph{krigging} methods
(see e.g. \cite{Wackern95}, \cite{cressie}) or spatial autoregressive
model (see for example \cite{cressie}) since it does not requires
any underlying model. It only requires the knowledge of the number
of the neighbors. We will see that the property of the \emph{inverse
regression} method provides a way of estimating this number.

\vskip .1in

This paper falls into the following parts. Section \ref{s:frame}
provides some notations and assumptions on the spatial process, as
well as some preliminar results on U-statistics. The estimation method
and the consistency results are presented in Section \ref{s:estim}.
Section \ref{s:fore} uses this estimate to forecast a spatial process.
Section \ref{sec:Conclus} is devoted to Conclusion. Proofs and the
technical lemmas are gathered in Section \ref{s:proof}.

\section{General setting and preliminary Results\label{s:frame}}

\subsection{Notations and assumptions}

Throughout all the paper, we will use the following notations.\\
 For all $b\in\mathbb{R}^{d}$ , $b^{(j)}$ will denote the $j^{th}$
component of the vector $b$;\\
 a point in bold $\mathbf{i}=(i_{1},...,i_{N})\in\mathbf{n}\in(\mathbb{N}^{*})^{N}$
will be referred to as a site, we will set $\mathbf{1}_{N}=(\underbrace{1,...,1}_{N\,\mbox{times}})$;
if $\mathbf{n}=(n_{1},...,n_{N})$, we will set $\widehat{\mathbf{n}}=n_{1}\times...\times n_{N}$
and write $\mathbf{n\rightarrow+\infty}$ if $\min_{i=1,...,N}n_{i}\,\,\mathbf{\rightarrow+\infty}$
and $\frac{n_{i}}{n_{k}}<C$ for some constant $C>0$.\\
 The symbol $\left\Vert .\right\Vert $ will denote any norm over
$\mathbb{R}^{d}$ , $\left\Vert u\right\Vert _{\infty}=\sup_{x}|u(x)|$
for some function $u$ and $C$ an arbitrary positive constant. If
$A$ is a set, let $1_{A}(x)=\left\{ \begin{array}{c}
1\,\,\textrm{ if }x\in A\\
0\,\,\textrm{ otherwise}\end{array}\right.$.\\
 The notation $W_{\mathbf{n}}=\mathcal{O}_{p}(V_{\mathbf{n}})$ (respectively
$W_{\mathbf{n}}=\mathcal{O}_{a.s}(V_{\mathbf{n}})$) means that $W_{\mathbf{n}}=V_{\mathbf{n}}S_{\mathbf{n}}$
for a sequence $S_{\mathbf{n}}$, which is bounded in \emph{probability}
(respectively \emph{almost surely}). \vskip .1in

We are interested in some $\mathbb{R}^{d}\times\mathbb{R}$-valued
stationary and measurable random field $Z_{\mathbf{i}}=(X_{\mathbf{i}},Y_{\mathbf{i}})$,
$\mathbf{i}\in(\mathbb{N}^{*})^{N}$, $(N,\, d\geq1)$ defined on
a probability space $(\Omega,\,\mathcal{A},\mathbf{P})$. Without
loss of generality, we consider estimations based on observations
of the process $(Z_{\mathbf{i}},\,\mathbf{i}\in\mathbb{Z}^{N})$ on
some rectangular set $\mathcal{I}_{\mathbf{n}}=\left\{ \mathbf{i}=(i_{1},...,i_{N})\in\mathbb{Z}^{N},\,1\leq i_{k}\leq n_{k},\, k=1,...,N\right\} $
for all $\mathbf{n}\in(\mathbb{N}^{*})^{N}$.\\
 Assume that the $Z_{\mathbf{i}}$'s have the same distribution as
$(X,Y)$ which is such that:

\begin{itemize}
\item the variable $Y$ has a density $f$. 
\item $\forall j=1,...,d$ each component $X^{(j)}$ of $X$, is such that
the pair $(X^{(j)},Y)$ admits an unknown density $f_{X^{(j)},Y}$
with respect to Lebesgue measure $\lambda$ over $\mathbb{R}^{2}$
and each $X^{(j)}$ is integrable. 
\end{itemize}

\subsection{Spatial dependency}

~~\\

As mentionned above, our model as any spatial model must take into
account spatial dependence between values at differents locations.
Of course, we could consider that there is a global linear relationships
between locations as it is generally done in spatial linear modeling,
we prefer to use a nonlinear spatial dependency measure. Actually,
in many circumstances the spatial dependency is not necessarly linear
(see \cite{arbialaf05}). It is, for example, the classical case where
one deals with the spatial pattern of extreme events such as in the
economic analysis of poverty, in the environmental science,... Then,
it is more appropriate to use a nonlinear spatial dependency measure
such as positive dependency (see \cite{arbialaf05}) or strong mixing
coefficients concept (see Tran \cite{tran}). In our case, we will
measure the spatial dependency of the concerned process by means of
\emph{$\alpha-$mixin}g and \emph{local dependency measure}.

\subsubsection{Mixing condition : }

$ $\\

The field $(Z_{\mathbf{i}})$ is said to satisfy a \emph{mixing condition}
if:

\begin{itemize}
\item there exists a function $\mathcal{X}:\mathbb{R}^{+}\to\mathbb{R}^{+}$
with $\mathcal{X}(t)\downarrow0$ as $t\to\infty$, such that whenever
$S,\, S'\subset(\mathbb{N}^{*})^{N}$ , 
\begin{align}
\alpha(\mathcal{B}(S),\mathcal{B}(S'))= & \sup_{A\in\mathcal{B}(S),\: B\in\mathcal{B}(S')}|P(B\cap C)-P(B)P(C)|\nonumber \\
 & \leq\psi(\mbox{Card}S,\,\mbox{Card}S')\,\mathcal{X}(\mbox{dist}(S,S'))\label{eq:mixingc}\end{align}
 where $\mathcal{B}(S)$(\emph{resp.} $\mathcal{B}(S')$) denotes
the Borel $\sigma-$fields generated by $(Z_{\mathbf{i}},\mathbf{i}\in S)$
(\emph{resp.} $(Z_{\mathbf{i}},\mathbf{i}\in S')$), $\mbox{Card}\, S$
(\emph{resp.} $\mbox{Card}S'$) the cardinality of $S$(\emph{resp.}
$S'$), $\mbox{dist}(S,S')$ the Euclidean distance between $S$ and
$S'$, and $\psi:\mathbb{N}^{2}\to\mathbb{R}^{+}$ is a symmetric
positive function nondecreasing in each variable. If $\psi\equiv1$,
then $Z_{\mathbf{i}}$ is called strong mixing. It is this latter
case which will be tackled in this paper and for all $v\geq0$, we
have \\
 \[
\alpha\left(v\right)=\sup_{\mathbf{i,\, j}\in\Bbb{R}^{N},\left\Vert \mathbf{i-}\mathbf{j}\right\Vert =v}\alpha\left(\sigma\left(Z_{\mathbf{i}}\right),\,\sigma\left(Z_{\mathbf{j}}\right)\right)\leq\mathcal{X}(v).\]

\item The process is said to be \emph{Geometrically Strong Mixing} (GSM)
if there exists a non-negative constant $\rho\in[0,1[$ such that
for all $u>0$, $\alpha(u)\leq C\rho^{u}$ . 
\end{itemize}
\begin{rem*}
A lot of published results have shown that the mixing condition (\ref{eq:mixingc})
is satisfied by many time series and spatial random processes (see
e.g. Tran \cite{tran}, Guyon \cite{guyon}, Rosenblatt \cite{Roseblat85},
Doukhan \cite{doukhan}). Moreover, the results presented in this
paper could be extended under additional technical assumptions to
the case, often considered in the literature, where $\psi$ satisfies:
\[
\psi(\mathbf{i},\,\mathbf{j})\leq c\,\min(\mathbf{i},\,\mathbf{j}),\,\,\,\,\,\forall\,\mathbf{i},\,\mathbf{j}\,\in\mathbb{N},\]
 for some constant $c>0$.

In the following, we will consider the case where $\alpha(u)\le Cu^{-\theta}$,
for some $\theta>0$. But, the results can be easly extend to the
\emph{GSM} case.
\end{rem*}
\bigskip{}

\subsubsection{Local dependency measure.}

$ $\\

In order to obtain the same rate of convergence as in the \emph{i.i.d}
case, one requires an other dependency measure, called a \emph{local
dependency measure}. Assume that

\begin{itemize}
\item For $\ell=1,...,d$, there exits a constant $\Delta>0$ such that
the pairs $(X_{\mathbf{i}}^{(\ell)},\, X_{\mathbf{j}})$ and ($(X_{\mathbf{i}}^{(\ell)},\, Y_{\mathbf{i}}),\,(X_{\mathbf{j}}^{(\ell)},\, Y_{\mathbf{j}})$)
admit densities $f_{\mathbf{i},\mathbf{j}}$ and $g_{\mathbf{i},\mathbf{j}}$,
as soon as $\mbox{dist}(\mathbf{i},\,\mathbf{j})>\Delta$, such that
\[
|f_{\mathbf{i},\,\mathbf{j}}\left(x,y\right)-f\left(x\right)f\left(y\right)|\leq C,\,\,\,\forall x,y\,\in\mathbb{R}\]
 \[
|g_{\mathbf{i},\,\mathbf{j}}\left(u,v\right)-g\left(u\right)g\left(v\right)|\leq C,\,\,\,\forall u,v\,\in\mathbb{R}^{2}\]
 for some constant $C\geq0$. 
\end{itemize}
\begin{rem*}
The link between the two dependency measures can be found in Bosq
\cite{bosq98b}. 

Note that if the second measure (as is name point out) is used to
control the local dependence, the first one is a kind of {}``asymptotic
dependency'' control.
\end{rem*}

\subsection{Results on U-statistics}

$ $\\

Let $(X_{n},\, n\geq1)$ be a sequence of real-valued random variables
with the same distribution as $F$. Let the functional: \[
\Theta(F)=\int_{\mathbb{R}^{m}}\, h(x_{1},x_{2},...,x_{m})dF(x_{1})...dF(x_{m}),\]
 where $m\in\mathbb{N}$, $h(.)$ is some measurable function, called
the kernel and $F$ is a distribution function from some given set
of distribution function. Without loss of generality, we can assume
that $h(.)$ is invariable by permutation. Otherwise, the transformation
$\frac{1}{m!}\sum_{1\leq i_{1}\neq i_{2}\neq...\neq i_{m}\leq n}h(x_{i_{1}},...,x_{i_{m}})$
will provide a symmetric kernel.

A $U-$statistic with kernel $h(.)$ of degree $m$ based on the sample
$(X_{i},\,1\leq i\leq n)$ is a statistic defined by: \begin{eqnarray*}
U_{\mathbf{n}} & = & \frac{(n-m)!}{n!}\sum_{1\leq i_{1}\neq i_{2}\neq...\neq i_{m}\leq n}h(X_{i_{1}},...,X_{i_{m}})\end{eqnarray*}
 It is said to be an $m-$order $U-$statistic. Let $h_{1}(x_{1})=\int_{\mathbb{R}^{m-1}}\, h(x_{1},x_{2},...,x_{m})\prod_{j=2}^{m}dF(x_{j}).$

The next Lemma is a consequence of Lemma 2.6 of Sun \& Chian \cite{Sun-Chian}.

\begin{lem}
\label{lem:Ustat-pro} ~~~Let $(X_{n},\, n\geq1)$ be a stationary
sequence of strongly mixing random variables. If there exists a positive
number $\delta$ and $\delta'$ $(0<\delta'<\delta)$ verifying $\gamma=\frac{6(\delta-\delta')}{(4+\delta)(2+\delta')}>1$
such that \begin{equation}
||h(X_{1},...,X_{m})||_{4+\delta}<\infty,\label{eq:Ustat-pro-Cond1}\end{equation}
 \begin{equation}
\int_{\mathbb{R}^{m}}\,|h(x_{1},...,x_{m})|^{4+\delta}\prod_{j=1}^{m}dF(x_{j})<\infty,\label{eq:Ustat-pro-Cond2}\end{equation}
 and $\alpha(n)=\mathcal{O}(n^{-3(4+\delta')/(2+\delta')})$ . Then,
\[
U_{n}=\Theta(F)+\frac{2}{n}\sum_{i=1}^{n}\,\left(h_{1}(X_{i})-\Theta(F)\right)+\mathcal{O}_{p}(\frac{1}{n}).\]

\end{lem}
~~~

To give strong consistency results, we need the following law of the
iterated logarithm of U-statistics:

\begin{lem}
\label{lem:Ustat-as} (Sun \& Chian, \cite{Sun-Chian}) Under the
same conditions of the previous lemma, we have \[
U_{n}-\Theta(F)=\frac{2}{n}\sum_{i=1}^{n}\,\left(h_{1}(X_{i})-\Theta(F)\right)+\mathcal{O}_{a.s}\left(\sqrt{\frac{\log\log n}{n}}\right).\]

\end{lem}
\begin{rem}
~~~ 

\begin{itemize}
\item In the following, we are dealing with a kernel $h(.)=K(\frac{.}{h_{\mathbf{n}}})$
which depends on $\mathbf{n}$. Actually, it is a classical approach
to use $U-$statistics result to get some assymptotic results of kernel
estimators, in the \emph{i.i.d} case, we refer for example Härdle
and Stoker \cite{Hardle-stoker89}. In fact, the dependence of $h_{\mathbf{n}}$
on $\mathbf{n}$ does not influence the asymptotical results presented
here. 
\end{itemize}
\end{rem}

\section{Estimation of the covariance of Inverse Regression Estimator\label{s:estim}}

We suppose that one deals with a random field $(Z_{\mathbf{i}},\,\mathbf{i}\in\mathbb{Z}^{N})$
which, corresponds, in the spatial regression case, to observations
of the form $Z_{\mathbf{i}}=(X_{\mathbf{i}},Y_{\mathbf{i}})$, $\mathbf{i}\in\mathbb{Z}^{N}$,
$(N\geq1)$ at different locations of a subset of $\mathbb{R}^{N}$,
$N\geq1$ with some dependency structure. Here, we are particularly
interested with the case where the locations take place in lattices
of $\mathbb{R}^{N}$. The general continuous case will be the subject
of a forthcoming work. \vskip .1in We deal with the estimation of
the matrix $\Sigma_{e}=\mathbf{var}\mathbf{E}(X|Y)$ based on the
observations of the process: $(Z_{\mathbf{i}},\,\mathbf{i}\in\mathcal{I}_{\mathbf{n}})$
; $\mathbf{n}\in\,\left(\mathbb{N}^{*}\right)^{N}$. In order to ensure
the existence of the matrix $\Sigma=\mathbf{var}\, X$ and $\Sigma_{e}=\mathbf{var}\,\mathbf{E}(X|Y)$,
we assume that $\mathbf{E}||X||^{4}<\infty$. For sake of simplicity
we will consider centered process so $\mathbf{E}X=0$.

To estimate model~(\ref{eq:mod}), as previously mentioned, one needs
to estimate the matrix $\Sigma^{-1}\Sigma_{e}$. On the one hand,
we can estimate the variance matrix $\Sigma$ by the empirical spatial
estimator, whose consistency will be easily obtained. On the other
hand, the estimation of the matrix $\Sigma_{e}$ is delicate since
it requires the study of the consistency of a suitable estimator of
the (inverse) regression function of $X$ given $Y$: \[
r(y)=\left\{ \begin{array}{cc}
\frac{\varphi(y)}{f(y)} & \textrm{if}\, f(y)\neq0;\\
\textrm{E}Y & \textrm{if}\, f(y)=0\end{array}\right.\mbox{where }\varphi(y)=\left(\int_{\mathbb{R}}x^{(i)}f_{X^{(i)},Y}(x^{(i)},y)dx,\,1\leq i\leq d\right),y\in\mathbb{R}.\]

An estimator of the \emph{inverse regression} function $r(.)$, based
on $(Z_{\mathbf{i}},\,\mathbf{i}\in\mathcal{I}_{\mathbf{n}})$ is
given by

\[
r_{\mathbf{n}}(y)=\left\{ \begin{array}{cc}
\frac{\varphi_{\mathbf{n}}(y)}{\, f_{\mathbf{n}}(y)} & \textrm{if}\, f_{\mathbf{n}}(y)\neq0,\\
\frac{1}{\hat{\mathbf{n}}}\sum_{\mathbf{i}\in\mathcal{I}_{\mathbf{n}}}Y_{\mathbf{i}} & \textrm{if}\, f_{\mathbf{n}}(y)=0,\end{array}\right.\]
 with for all $y\in\mathbb{R}$, \[
f_{\mathbf{n}}(y)=\frac{1}{\widehat{\mathbf{n}}h_{\mathbf{n}}}\sum_{\mathbf{i}\in\mathcal{I}_{\mathbf{n}}}K\left(\frac{y-Y_{\mathbf{i}}}{h_{\mathbf{n}}}\right)\]
 \[
\varphi_{\mathbf{n}}(y)=\frac{1}{\widehat{\mathbf{n}}h_{\mathbf{n}}}\sum_{\mathbf{i}\in\mathcal{I}_{\mathbf{n}}}X_{\mathbf{i}}K\left(\frac{y-Y_{\mathbf{i}}}{h_{\mathbf{n}}}\right),\]
 where $f_{\mathbf{n}}$ is a kernel estimator of the density, $K:\,\mathbb{R}^{d}\rightarrow\mathbb{R}$
is a bounded integrable kernel such that $\int K\left(x\right)dx=1$
and the bandwidth $h_{\mathbf{n}}\geq0$ is such that $\lim_{n\rightarrow+\infty}h_{\mathbf{n}}=0$.

The consistency of the estimators $f_{\mathbf{n}}$ and $r_{\mathbf{n}}$
has been studied by Carbon et al \cite{carbon-tran-wu}. To prevent
small-valued density observations $y$, we consider the following
density estimator: \[
f_{e,\mathbf{n}}(y)=\max(e_{\mathbf{n}},f_{\mathbf{n}}(y))\]
 where $(e_{\mathbf{n}})$ is a real-valued sequence such that $\lim_{\mathbf{n}\to\infty}e_{\mathbf{n}}=0$.
Then, we consider the corresponding estimator of $r$ \[
r_{e,\mathbf{n}}(y)=\frac{\varphi_{\mathbf{n}}(y)}{\, f_{e,\mathbf{n}}(y)}.\]

Finally, for $\overline{X}=\frac{1}{\hat{\mathbf{n}}}\sum_{\mathbf{i}\in\mathcal{I}_{\mathbf{n}}}X_{\mathbf{i}}$
we consider the estimator of $\Sigma_{e}$:

\[
\Sigma_{e,\mathbf{n}}=\frac{1}{\hat{\mathbf{n}}}\sum r_{e,\mathbf{n}}(Y_{\mathbf{i}})\, r_{e,\mathbf{n}}(Y_{\mathbf{i}})^{T}-\overline{X}\,\overline{X}^{T}.\]

We aim at proving the consistency of the empirical variance associated
to this estimator.

\begin{rem*}
~~ 
Here, we consider as estimator of the density $f$, $f_{e,\mathbf{n}}=\max(e_{\mathbf{n}},f_{\mathbf{n}})$,
to avoid small values. There are other alternatives such as $f_{e,\mathbf{n}}=f_{\mathbf{n}}+e_{\mathbf{n}}$
or $f_{e,\mathbf{n}}=\max\{(f_{\mathbf{n}}-e_{\mathbf{n}}),0\}$.

\end{rem*}
\vskip .1in

\subsection{Weak consistency}

~~ In the following, for a fixed $\eta>0$ and a random variable
$Z$ in $\mathbb{R}$$^{d}$, we will use the notation $\left\Vert Z\right\Vert _{\eta}=\mathbf{E}(||Z||^{\eta})^{1/\eta}$.

In this section, we will make the following technical assumptions

\begin{equation}
\left\Vert \frac{r(Y)}{f(Y)}\right\Vert _{4+\delta_{1}}<\infty,\,\,\,\mbox{for some }\delta_{1}>0\label{eq:Assump-rrT/f}\end{equation}
 and \begin{equation}
\left\Vert \frac{r(Y)\,}{f(Y)}\mathbf{1}_{\{f(Y)\leq e_{\mathbf{n}}\}}\right\Vert _{2}=\mathcal{O}\left(\frac{1}{\hat{\mathbf{n}}^{\frac{1+\delta}{2}}}\right).\,\,\,\,\mbox{for some}\,\,1>\delta>0.\end{equation}
These assumptions are the spatial counterparts of respectively $\left\Vert r(Y)\right\Vert _{4+\delta}<\infty$
and $\left\Vert r(Y)\,\mathbf{1}_{\{f(Y)\leq e_{\mathbf{n}}\}}\right\Vert _{2}=\mathcal{O}\left(\frac{1}{\hat{\mathbf{n}}^{\frac{1}{4}+\delta}}\right)$
needed in the \emph{i.i.d} case.

We also assume some regularity conditions on the functions: $K(.)$,
$f(.)$ and $r(.)$:

\begin{itemize}
\item The kernel function $K(.):\,\mathbb{R}\to\mathbb{R}^{+}$ is a $k-$order
kernel with compact support and satisfying a Lipschitz condition $\left|K\left(x\right)-K\left(y\right)\right|\leq C|x-y|$ 
\item $f(.)$ and $r(.)$ are functions of $C^{k}(\mathbb{R})$ ($k\geq2$)
such that $\sup_{y}|f^{(k)}(y)|<C_{1}$ and $\sup_{y}||\varphi^{(k)}(y)||<C_{2}$
for some constants $C_{1}$ and $C_{2}$, 
\end{itemize}
Set $\Psi_{\mathbf{n}}=h_{\mathbf{n}}^{k}+\frac{\sqrt{\log\hat{\mathbf{n}}}}{\sqrt{\hat{\mathbf{n}}h_{\mathbf{n}}}}\,.$

\begin{thm}
\label{thm:conv-prob} Assume that $\alpha(t)\le Ct^{-\theta}$, $t>0$,
$\theta>2N$ and $C>0$. If $E(||X||)<\infty$ and $\psi(.)=\mathbf{E}(||X||^{2}|Y=.)$
is continuous. Then for a choice of $h_{\mathbf{n}}$ such that $\widehat{\mathbf{n}}h_{\mathbf{n}}^{3}(\log\widehat{\mathbf{n}})^{-1}\rightarrow0$
and $\widehat{\mathbf{n}}h_{\mathbf{n}}^{\theta_{1}}(\log\widehat{\mathbf{n}})^{-1}\rightarrow\infty$
with $\theta_{1}=\frac{4N+\theta}{\theta-2N}$,~ then, we get \[
\Sigma_{e,\mathbf{n}}-\Sigma_{e}=\mathcal{O}_{p}\left(h_{\mathbf{n}}^{k}+\frac{\Psi_{\mathbf{n}}^{2}}{e_{\mathbf{n}}^{2}}\right)\]
 
\end{thm}
\begin{cor}
\label{cor:cov-prob-vitesse opt} Under Assumptions of Theorem \ref{thm:conv-prob}
with $h\simeq n^{-c_{1}}$ , $e_{n}\simeq n^{-c_{2}}$ for some positive
constants $c_{1}$ and $c_{2}$ such that $\frac{c_{2}}{k}+\frac{1}{4k}<c_{1}<\frac{1}{2}-2c_{2},$
we have \[
\Sigma_{e,\mathbf{n}}-\Sigma_{e}=o_{p}\left(\frac{1}{\sqrt{\hat{\mathbf{n}}}}\right).\]

\end{cor}
~

\begin{cor}
\label{thm:(Central-limit-theorem)}(Central limit theorem) Under
previous assumptions, we have \[
\sqrt{\hat{\mathbf{n}}}\,\left(\Sigma_{e,\mathbf{n}}-\Sigma_{e}\right)\stackrel{\mathcal{L}}{\rightarrow}\Lambda\]
 where $\Lambda$ is a zero-mean gaussian on the space of $d$-order
matrix with covariance $\mathbf{var}\left(r(Y)r(Y)^{T}\right).$ 
\end{cor}
~~

\subsection{Strong consistency}

~~~

Here we study the case where the response, $Y$ takes values in some
compact set. We replace the assumption $\left\Vert \frac{r(Y)\,}{f(Y)}\mathbf{1}_{\{f(Y)\leq e_{\mathbf{n}}\}}\right\Vert _{2}=\mathcal{O}\left(\frac{1}{\hat{\mathbf{n}}^{\frac{1}{2}+\delta}}\right)$
by $\mathbf{E}\left(\exp\left(\left\Vert r(Y)\right\Vert \,\mathbf{1}_{\{f(Y)\leq e_{\mathbf{n}}\}}\right)\right)=\mathcal{O}\left(\hat{\mathbf{n}}^{-\xi}\right)$
for some $\xi>0$. : $\mathbf{E}\exp\gamma||X||<\infty$ for some
constant $\gamma>0$.

\begin{thm}
\label{thm:Conv-a.s} If $(Z_{\mathbf{u}})$ is GSM, for a choice
of $h_{\mathbf{n}}$ such that $\widehat{\mathbf{n}}h_{\mathbf{n}}^{3}(\log\widehat{\mathbf{n}})^{-1}\rightarrow0$
and $\hat{\mathbf{n}}\, h_{\mathbf{n}}(\log\hat{\mathbf{n}})^{-2N-1}\to\infty$.
Assume also that $\inf_{S}f(y)>0$ for some compact set $S$, then
under the Assumptions of Lemma \ref{lem:Ustat-pro}, we have: \[
\Sigma_{e,\mathbf{n}}-\Sigma_{e}=\mathcal{O}_{a.s}\left(h_{\mathbf{n}}^{k}+\frac{\Psi_{\mathbf{n}}^{2}}{e_{\mathbf{n}}^{2}}\right)\,.\]

\end{thm}
\begin{cor}
\label{cor:Conv-a.s-opt-rate}Under previous Assumptions, with $h_{\mathbf{n}}\simeq\left(\hat{\mathbf{n}}\right)^{-c_{1}}$,
$e_{\mathbf{n}}\simeq\mathbf{\hat{n}}^{-c_{2}}$ for some positive
constants $c_{1}$ and $c_{2}$ such that $\frac{c_{2}}{k}+\frac{1}{4k}\leq c_{1}<\frac{1}{2}-2c_{2},$
we get \[
\Sigma_{e,\mathbf{n}}-\Sigma_{e}=o_{a.s}\left(\sqrt{\frac{\log\log\hat{\mathbf{n}}}{\hat{\mathbf{n}}}}\right).\]

\end{cor}
As mentionned previously, the eigenvectors associated with the positive
eigenvalues of $\Sigma_{\mathbf{n}}^{-1}\Sigma_{e,\,\mathbf{n}}$
provide an estimation of the EDR space. Classically, weak and strong
consistency results concerning the estimation of the EDR space are
obtained by using the previous consistency respectively of the $\Sigma$
and $\Sigma_{e}$ and the theory of perturbation as for example in
\cite{Fang-Zhu96}.

\section{Spatial inverse methode for spatial prediction\label{s:fore}}

\subsection{Prediction of a spatial process}

$ $\\
Let $(\xi_{\mathbf{n}},\,\mathbf{n}\in(\mathbb{N}^{*})^{N})$ be a
$\mathbb{R}-$valued strictly stationary random spatial process, assumed
to be observed over a subset $\mathcal{O}_{\mathbf{n}}\subset\mathcal{I}_{\mathbf{n}}$
($\mathcal{I}_{\mathbf{n}}$ is a rectangular region as previously
defined for some $\mathbf{n}\in(\mathbb{N}^{*})^{N}$). Our aim is
to predict the square integrable value, $\xi_{\mathbf{i}_{0}}$, at
a given site $\mathbf{i}_{0}\in I_{\mathbf{n}}-\mathcal{O}_{\mathbf{n}}$.
In practice, one expects that $\xi_{\mathbf{i}_{0}}$ only depends
on the values of the process on a bounded vicinity set (as small as
possible) $\mathcal{V}_{\mathbf{i}_{0}}\subset\mathcal{O}_{\mathbf{n}}$;
i.e that the process $(\xi_{\mathbf{i}})$ is (at least locally) a
Markov Random Field (MRF) according to some system of vicinity. Here,
we will assume (without loss of generality) that the set of vicinity
$(\mathcal{V}_{\mathbf{j}},\mathbf{j}\in(\mathbb{N}^{*})^{N})$ is
defined by $\mathcal{V}_{\mathbf{j}}$ of the form $\mathbf{j}+\mathcal{V}$
(call vicinity prediction in Biau and Cadre \cite{biau-cadre}). Then
it is well known that the minimum \emph{mean-square error of prediction}
of $\xi_{\mathbf{i}_{0}}$ given the data in $\mathcal{V}_{\mathbf{i}_{0}}$
is \[
E(\xi_{\mathbf{i}_{0}}|\xi_{\mathbf{i}},\mathbf{i}\in\mathcal{V}_{\mathbf{i}_{0}})\]
 and we can consider as predictor any $d-$dimensional vector (where
$d$ is the cardinal of $\mathcal{V}$) of elements of $\mathcal{V}_{\mathbf{i}_{0}}$
concatenated and ordered according to some order. Here, we choose
the vector of values of $(\xi_{\mathbf{n}})$ which correspond to
the $d-$nearest neighbors: for each $\mathbf{i}\in\mathbb{Z}^{N}$,
we consider that the predictor is the vector $\xi_{\mathbf{i}}^{d}=(\xi_{\mathbf{i}(k)};\,1\leq k\leq d)$
where $\mathbf{i}(k)$ is the $k-$th nearest neighbor of $\mathbf{i}$.
Then, our problem of prediction amounts to estimate :

\[
m(x)=E(\xi_{\mathbf{i}_{0}}|\xi_{\mathbf{i}_{0}}^{d}=x).\]
 For this purpose we construct the \emph{associated process}: \[
Z_{\mathbf{i}}=(X_{\mathbf{i}},Y_{\mathbf{i}})=(\xi_{\mathbf{i}}^{d},\xi_{\mathbf{i}}),\,\mathbf{i}\in\mathbb{Z}^{N}\]
 and we consider the estimation of $m(.)$ based on the data $(Z_{\mathbf{i}},\in\mathcal{O}_{\mathbf{n}})$
and the model (\ref{eq:mod}). Note that the linear approximation
of $m(.)$ leads to linear predictors. The available literature on
such spatial linear models (we invite the reader think of the \emph{kriging
method }or spatial auto-regressive method\emph{)} is relatively abundant,
see for example, Guyon \cite{guyon}, Anselin and Florax \cite{anselin},
Cressie \cite{cressie}, Wackernagel \cite{Wackern95}. In fact, the
linear predictor is the optimal predictor (in mimimun mean square
error meaning) when the random field under study is \emph{Gaussian}.
Then, linear techniques for spatial predicition, give unsatisfactory
results when the the process is not \emph{Gaussian}. In this latter
case, other approaches such as \emph{log-normal kriging} or the \emph{trans-Gaussian
kriging} have been introduced. These methods consist in transforming
the original data into a Gaussian distributed data. But, such methods
lead to outliers which appear as an effect of the heavy-tailed densities
of the data and cannot be delete. Therefore, a specific consideration
is needed. This can be done by using, for example, a nonparametric
model. That is what is proposed by Biau and Cadre \cite{biau-cadre}
where a predictor based on \emph{kernel methods} is developped. But,
This latter (the kernel nonparametric predictor) as all kernel estimator
is submitted to the so-called \emph{dimension curse} and then is penalized
by $d$ ($=\mbox{card}(\mathcal{V})$), as highlighted in Section
1. Classically, as in Section 1, one uses dimension reduction such
as the \emph{inverse regression} method, to overcome this problem.
We propose here an adaptation of the \emph{inverse regression} method
to get a \emph{dimension reduction predictor} based on model (\ref{eq:mod}):
\begin{equation}
\xi_{\mathbf{i}}=g(\Phi.\xi_{\mathbf{i}}^{d}).\label{eq:predmod}\end{equation}

\begin{rem}
~~~
\begin{enumerate}
\item To estimate this model, we need to check the SIR condition in the
context of prediction i.e: \emph{$X$ is such that for all vector
$b$ in $\mathbb{R}^{d}$, there exists a vector $B$ of $\mathbb{R}^{D}$
such that $\mathbf{E}(b^{T}X|\Phi.X)=B^{T}(\Phi.X)$, }is verify if
the process $(\xi_{\mathbf{i}})$ is a spatial elliptically distributed
process such as Gaussian random field\emph{.} 
\item In the time series forecasting problem, {}``inverse regression''
property can be an {}``handicap'', since then, one needs to estimate
the expectation of the \emph{{}``future''} given the \emph{{}``past''}.
So, the process under study must be reversible. The flexibility that
provide spatial modelling overcome this default since as mentioned
in the introduction, the notion of past, present and future does not
exist.
\end{enumerate}
\end{rem}
~~~

At this stage, one can use the method of estimation of the model (\ref{eq:mod})
given in Section 1 to get a predictor. Unfortunately (as usually in
prediction problem) $d$ is unknown in practice. So, we propose to
estimate $d$ by using the fact that we are dealing both with a Markov
property and \emph{inverse regression} as follows.

\subsection{Estimation of the number of neighbors necessary for prediction}

$ $\\

Note that we suppose that the underline process is a stationary Markov
process with respect to the $d-$neighbors system of neighborhood,
so the variables $\xi_{\mathbf{i}(k)}$ and $\xi_{\mathbf{i}}$ are
independent as soon as $k>d$ and \[
\mathbf{E}(\xi_{\mathbf{i}(k)}|\xi_{\mathbf{i}}=y)=0\]
 (since $(\xi_{\mathbf{i}})$ is a stationary zero mean process).

Futhermore since our estimator (of model (\ref{eq:mod})) is based
on estimation of $\mathbf{E}(X|Y=y)=\mathbf{E}(\xi_{\mathbf{i}}^{d}|\xi_{\mathbf{i}}=y)=(\mathbf{E}(\xi_{\mathbf{i}(k)}|\xi_{\mathbf{i}}=y);\,1\leq k\leq d)$,
that allows us to keep only the neighbors $\xi_{\mathbf{i}(k)}$ for
which $\mathbf{E}(\xi_{\mathbf{i}(k)}|\xi_{\mathbf{i}}=y)\neq0$.
Then, an estimation of $d$ is obtained by estimation of $\mbox{argmin}_{k}\mathbf{E}(\xi_{\mathbf{i}(k)}|\xi_{\mathbf{i}}=y)=0$.
We propose the following algorithm to get this estimator.

\subsection*{Algorithm for estimation of $d$, the number of neighbors.}

~~~~

\begin{enumerate}
\item Initialization: specify a parameter $\delta>0$ (small) and fix a
site $\mathbf{j}_{0}$; set $k=1$. 
\item compute $r_{\mathbf{n}}^{(k)}(y)=\frac{{\displaystyle \sum_{\mathbf{i}\in\mathcal{O}_{\mathbf{n}},\mathcal{V}_{\mathbf{j}_{0}}\subset\mathcal{O}_{\mathbf{n}}}}\xi_{\mathbf{i}(k)}\, K_{h_{\mathbf{n}}}\left(y-\xi_{\mathbf{i}}\right)}{{\displaystyle \sum_{\mathbf{i}\in\mathcal{O}_{\mathbf{n}},\mathcal{V}_{\mathbf{j}_{0}}\subset\mathcal{O}_{\mathbf{n}}}}K_{h_{\mathbf{n}}}\left(y-\xi_{\mathbf{i}}\right)}$,
the kernel estimate of $r^{(k)}(y)=\mathbf{E}(X^{(k)}|Y=y)$ 
\item if $|(r_{\mathbf{n}}^{(k)}(y)|>\delta$, then $k=k+1$ and continue
with Step 2; otherwise terminate and $d=k$. 
\end{enumerate}
Then, we can compute a predictor based on $d=k$:

\subsection{The dimension reduction predictor}

$ $\\

To get the predictor, we suggest the following algorithm:

\begin{enumerate}
\item compute \[
r_{\mathbf{n}}^{*}(y)=\frac{{\displaystyle \sum_{\mathbf{i}\in\mathcal{O}_{\mathbf{n}},\mathcal{V}_{\mathbf{i}_{0}}\subset\mathcal{O}_{\mathbf{n}}}}\xi_{\mathbf{i}}^{d}\, K_{h_{\mathbf{n}}}\left(y-\xi_{\mathbf{i}}\right)}{{\displaystyle \sum_{\mathbf{i}\in\mathcal{O}_{\mathbf{n}},\mathcal{V}_{\mathbf{i}_{0}}\subset\mathcal{O}_{\mathbf{n}}}}K_{h_{\mathbf{n}}}\left(y-\xi_{\mathbf{i}}\right)}\]

\item compute\[
\Sigma_{e,\mathbf{n}}=\frac{1}{\hat{\mathbf{n}}}{\displaystyle \sum_{\mathbf{i}\in\mathcal{O}_{\mathbf{n}},\mathcal{V}_{\mathbf{i}_{0}}\subset\mathcal{O}_{\mathbf{n}}}}r_{e,\mathbf{n}}^{*}(Y_{\mathbf{i}})\, r_{e,\mathbf{n}}^{*}(Y_{\mathbf{i}})^{T}-\overline{X}\,\overline{X}^{T}.\]

\item Do the \emph{principal component analisys} of $\Sigma_{\mathbf{n}}^{-1}\Sigma_{e,\,\mathbf{n}}$
both to get a basis of $\mbox{Im}(\Sigma_{\mathbf{n}}^{-1}\Sigma_{e,\,\mathbf{n}})$
and estimation of the $D$, the dimension of $\mbox{Im}(\Phi)$ as
suggested in the next remark 
\item compute the predictor: \[
\hat{\xi}_{\mathbf{i}_{0}}=g_{\mathbf{n}}^{*}(\Phi_{\mathbf{n}}^{*}.X_{\mathbf{i}_{0}}).\]

based on data $(Z_{\mathbf{i}},\mathbf{i}\in\mathcal{O}_{\mathbf{n}})$;
where $g_{\mathbf{n}}^{*}$ is the kernel estimate: \[
g_{\mathbf{n}}^{*}(x)=\frac{{\displaystyle \sum_{\mathbf{i}\in\mathcal{O}_{\mathbf{n}},\mathcal{V}_{\mathbf{i}_{0}}\subset\mathcal{O}_{\mathbf{n}}}}\xi_{\mathbf{i}}\, K_{h_{\mathbf{n}}}\left(\Phi_{\mathbf{n}}^{*}(x-\xi_{\mathbf{i}}^{d})\right)}{{\displaystyle \sum_{\mathbf{i}\in\mathcal{O}_{\mathbf{n}},\mathcal{V}_{\mathbf{i}_{0}}\subset\mathcal{O}_{\mathbf{n}}}}K_{h_{\mathbf{n}}}\left(\Phi_{\mathbf{n}}^{*}(x-\xi_{\mathbf{i}}^{d})\right)}\forall x\in\mathbb{R}^{d}.\]

\end{enumerate}
\medskip{}

\begin{rem}
~~ 
\begin{enumerate}
\item The problem of estimation of $D$ in step (4) is a classical problem
in dimension reduction problems. Several ways exist in the literature.
One can for example use the eigenvalues representation of the matrix
$\Sigma_{\mathbf{n}}^{-1}\Sigma_{e,\,\mathbf{n}}$, the measure of
distance between spaces as in Li \cite{Li} or the selection rule
of Ferr\'e \cite{Ferre98}.
\item Consitency on the convergence of $\hat{\xi}_{\mathbf{i}_{0}}$ to
$\xi_{\mathbf{i}_{0}}$ can be obtained by sketching both result of
Section \ref{s:estim} and results Biau and Cadre \cite{biau-cadre}. 
\end{enumerate}

\section{Conclusion\label{sec:Conclus} }

In this work, we have proposed two dimension reduction methods for
spatial modeling. The first one is a dimension reduction for spatial
regression. It is a natural extension of the idea of Li \cite{Li}
(called Inverse Regression method) for spatially dependent variables
under strong mixing condition. Then, on one hand, we can say that
is a good alternative to spatial linear regression model since the
link between the variables $X$ and $Y$ is not necessarly linear.
Futhermore, as raises Li \cite{Li}, any linear model can be seen
as a particular case of model (\ref{eq:mod}) with $g$ being the
identity function and $D=1$. On the other hand, as in the \emph{i.i.d}.
case, it requieres less data for calculus than spatial non-parametric
regression methods.

The second method that we have studied here deals with spatial prediction
modelling. Indeed, it is more general than \emph{kriging method} were
the gaussian assumption on the $X$ is needed. Here, we requier that
$X$ belongs to a larger class of random variables (that obey to Li
\cite{Li}'s condition recalled in the introduction). Futhermore,
our spatial prediction method has the ease of implementation property
of the \emph{inverse regression} methods. Then, for example, it allows
to estimate the number of neighbors need to predict. That cannot do
the non-parametric prediction method of Biau and Cadre \cite{biau-cadre}. 

We have presented here the theoretical framework of our techniques.
The next step is to apply them on real data. It is the subjet of works
under development.
\end{rem}

\section{Proofs and Technical Results\label{s:proof}}

\subsection{Deviation Bounds \label{s:prem} }

To show the strong consistency results, we will use the following
Bernstein type deviation inequality:

\begin{lem}
\label{lem:Dabo-Yao06} Let (\emph{$\zeta_{\mathbf{v}},\,\mathbf{v}\in\mathbb{N}^{N})$}
be a zero-mean real-valued random spatial process such that each $\mathbf{v}\in(\mathbb{N}^{*})^{N}$
there exists $c>0$ verifying \begin{equation}
\mathbf{E}\,|\zeta_{\mathbf{v}}|^{k}\leq k!\, c^{k-2}\,\mathbf{E}\,|\zeta_{\mathbf{v}}|^{2}\,,\,\,\forall\, k\geq2\label{eq:hyplemdabo-yao}\end{equation}
 for some constant $c>0$. Let $S_{\mathbf{n}}=\sum_{\mathbf{v}\in I_{\mathbf{n}}}\zeta_{\mathbf{v}}$.
Then for each $r\in[1,+\infty]$ and each $\mathbf{n}\in(\mathbb{N}^{*})^{N}$and
$\mathbf{q}\in(\mathbb{N}^{*})^{N}$ such that $1\leq q_{i}\leq\frac{n_{i}}{2}$
and each $\varepsilon>0$,\begin{equation}
P(|S_{\mathbf{n}}|>\hat{\mathbf{n}}\varepsilon)\leq2^{N+1}\textrm{exp}\left(-\frac{\hat{\mathbf{q}}\varepsilon^{2}}{4(\, M_{2}^{2}+2^{N}c\varepsilon)}\right)+2^{N}\times\hat{\mathbf{q}}\times11\left(1+\frac{4c\, p^{N}M_{2}^{2/r}}{\varepsilon}\right)^{\frac{r}{2r+1}}\alpha([p])^{2r/(2r+1)}\label{eq:large-deviation-inequal}\end{equation}
 where $M_{2}^{2}=\sup_{\mathbf{v}\in\mathcal{I}_{\mathbf{n}}}\mathbf{E}\zeta_{\mathbf{v}}^{2}$. 
\end{lem}
\begin{rem}
Actually, this result is an extension of Lemma 3.2 of Dabo-Niang and
Yao \cite{dabo-yao07} for bounded processes. This extension is necessary
since in the problem of our interessed, assuming the boundness of
the processes amounts to assume that the $X_{\mathbf{i}}$'s are bounded.
It is a restrictive condition which (generally) is incompatible with
the cornerstone condition of the \emph{inverse regression} (if $X$
is elliptically distributed for example). 
\end{rem}
We will use the following lemma to get the weak consistency and a
law of iterated \textit{\emph{of the logarithm as well as for the
matrix $\Sigma$ (as we will see immediately) than for the matrix
$\Sigma_{e}$ (see the proofs of results of Section \ref{s:estim}).}}

\begin{lem}
\label{lem: gLaw-of-large} Let $\{X_{\mathbf{n}},\,\mathbf{n}\in\mathbb{N}^{N}\}$
be a zero-mean stationary spatial process sequence, of strong mixing
random variables. 
\end{lem}
\begin{enumerate}
\item If $\mathbf{E}||X||^{2+\delta}<+\infty$ and $\sum\alpha(\hat{\mathbf{n}})^{\frac{\delta}{2+\delta}}<\infty$,
for some $\delta>0$. Then,\[
\frac{1}{\hat{\mathbf{n}}}\sum_{\mathbf{i}\in\mathcal{I}_{\mathbf{n}}}X_{\mathbf{i}}=\mathcal{O}_{p}\left(\frac{1}{\hat{\mathbf{n}}}\right).\]

\item If $\mathbf{E}||X||^{2+\delta}<+\infty$ and $\sum\alpha(\hat{\mathbf{n}})^{\frac{\delta}{2+\delta}}<\infty$,
for some $\delta>0$. Then,\[
\sqrt{\hat{\mathbf{n}}}(\frac{1}{\hat{\mathbf{n}}}\sum_{\mathbf{i}\in\mathcal{I}_{\mathbf{n}}}X_{\mathbf{i}})/\sigma\to\mathcal{N}(0,1)\]
 with $\sigma^{2}=\sum_{i\in\mathbb{Z}^{N}}\mbox{cov}(X_{k},X_{\mathbf{i}})$ 
\item If $\mathbf{E}\exp\gamma||X||<\infty$ for some constant $\gamma>0$,
if for all $u>0$ , $\alpha(u)\leq a\rho^{-u}$ , $0<\rho<1$ or $\alpha(u)=C.u^{-\theta}$,
$\theta>N$ then,\[
\frac{1}{\hat{\mathbf{n}}}\sum_{\mathbf{i}\in\mathcal{I}_{\mathbf{n}}}X_{\mathbf{i}}=o_{a.s}\left(\sqrt{\frac{\log\log\hat{\mathbf{n}}}{\hat{\mathbf{n}}}}\right).\]

\end{enumerate}
\begin{rem}
~~ 
\begin{itemize}
\item \textit{The first result is obtained by using covariance inequality
for strong mixing processes (see Bosq \cite{bosq98b}). Actually,
it suffices to enumerate the $X_{i}$ 's into an arbitrary order and
sketch the proof in Theorem 1.5 of Bosq \cite{bosq98b}.} 
\item \textit{The law of the iterated of the logarithm holds by applying
the previous Lemma \ref{lem:Dabo-Yao06} with $\varepsilon=\eta\,\sqrt{\frac{\log\log\hat{\mathbf{n}}}{\mathbf{\hat{\mathbf{n}}}}}$,
$\eta>0$ and $\hat{q}=\left[\,\frac{\hat{\mathbf{n}}}{\log\log\hat{\mathbf{n}}}\right]+1$.} 
\end{itemize}
\end{rem}
\vskip .1in

\subsection{Consistency of the inverse regression}

In Section \ref{s:estim}, we have seen that the results are based
on consistency results of the function $r(.)$ which are presented
now under some regularity conditions on the functions: $K(.)$, $f(.)$
and $r(.)$.

\begin{itemize}
\item The kernel function $K(.):\,\mathbb{R}\to\mathbb{R}^{+}$ is a $k-$order
kernel with compact support and satisfying a Lipschitz condition $\left|K\left(x\right)-K\left(y\right)\right|\leq C|x-y|$ 
\item $f(.)$ and $r(.)$ are functions of $C^{k}(\mathbb{R})$ ($k\geq2$)
such that $\sup_{y}|f^{(k)}(y)|<C_{1}$ and $\sup_{y}||\varphi^{(k)}(y)||<C_{2}$
for some constants $C_{1}$ and $C_{2}$, 
\end{itemize}
we have convergence result:

\begin{lem}
\label{thm:Carbon-et-alPro} Suppose $\alpha(t)\le Ct^{-\theta}$,
$t>0$, $\theta>2N$ and $C>0$. If $\widehat{\mathbf{n}}h_{\mathbf{n}}^{3}(\log\widehat{\mathbf{n}})^{-1}\rightarrow0$,
$\widehat{\mathbf{n}}h_{\mathbf{n}}^{\theta_{1}}(\log\widehat{\mathbf{n}})^{-1}\rightarrow\infty$
with $\theta_{1}=\frac{4N+\theta}{\theta-2N}$,~ then 
\end{lem}
\begin{enumerate}
\item (see, \cite{carbon-tran-wu}) \begin{equation}
\mbox{sup}_{y\in\mathbb{R}}|f_{\mathbf{n}}(y)-f(y)|=\mathcal{O}_{p}\left(\Psi_{\mathbf{n}}\,\right).\label{eq:Carbonfp}\end{equation}

\item Furthermore, if $E(||X||)<\infty$ and $\psi(.)=\mathbf{E}(||X||^{2}|Y=.)$
is continuous, then 
\end{enumerate}
\begin{equation}
\mbox{sup}_{y\in\mathbb{R}}||\varphi_{\mathbf{n}}(y)-\varphi(y)||=\mathcal{O}_{p}\left(\Psi_{\mathbf{n}}\right).\label{eq:Carbonphip}\end{equation}

~~~

\begin{rem}
Actually, only the result (\ref{eq:Carbonfp}) is shown in Carbon
et al \cite{carbon-tran-wu} but the result (\ref{eq:Carbonphip})
is easily obtained by noting that for all $\varepsilon>0$,\[
\mathbf{P}(\mbox{sup}_{y\in\mathbb{R}}||\varphi_{\mathbf{n}}(y)-\mathbf{E}\varphi_{\mathbf{n}}(y)||>\varepsilon)\leq\frac{\mathbf{E}||X||}{a_{\mathbf{n}}}+\mathbf{P}(\mbox{sup}_{y\in\mathbb{R}}||\varphi_{\mathbf{n}}(y)-\mathbf{E}\varphi_{\mathbf{n}}(y)||>\varepsilon,\,\forall i,\,||X_{i}||\leq a_{\mathbf{n}})\]
 with $a_{\mathbf{n}}=\eta\,(\log\hat{\mathbf{n}})^{1/4},\,\eta>0$. 
\end{rem}
\begin{lem}
\label{thm:Carbon-et-alAS} If $(Z_{\mathbf{u}})$ is GSM, $\widehat{\mathbf{n}}h_{\mathbf{n}}^{3}(\log\widehat{\mathbf{n}})^{-1}\rightarrow0$
and $\hat{\mathbf{n}}\, h_{\mathbf{n}}(\log\hat{\mathbf{n}})^{-2N-1}\to\infty,$
then\begin{equation}
\mbox{sup}_{y\in\mathbb{R}}|f_{\mathbf{n}}(y)-f(y)|=\mathcal{O}_{a.s}\left(\Psi_{\mathbf{n}}\right).\label{eq:Carbonfa.s}\end{equation}

Furthermore, if $\mathbf{E}\,(\exp\gamma\,||X||)<\infty$ for some
$\gamma>0$ and $\psi(.)=\mathbf{E}(||X||^{2}|Y=.)$ is continuous,
then \begin{equation}
\mbox{sup}_{y\in\mathbb{R}}||\varphi_{\mathbf{n}}(y)-\varphi(y)||=\mathcal{O}_{a.s}\left(\Psi_{\mathbf{n}}\right).\label{eq:Carbonphia.S}\end{equation}

\end{lem}
\begin{rem*}
The equality (\ref{eq:Carbonfa.s}) is due to Carbon et al \cite{carbon-tran-wu}.
The proof of the equality (\ref{eq:Carbonphia.S}) is obtained applying
Lemma \ref{lem:Dabo-Yao06} and sketching the proofs of Theorem 3.1
and 3.3 of Carbon et al \cite{carbon-tran-wu}. Then it is omitted. 
\end{rem*}
We will need the following lemma and the spatial block decomposition:

\begin{lem}
\label{lem:(-Bosq,-1997)}(Bradley's Lemma in Bosq \cite{bosq97})

Let $(X,Y)$ be an $\mathbb{R}^{d}\times\mathbb{R}-$valued random
vector such that $Y\in\mathbf{L}^{r}(P)$ for some $r\in[1,+\infty]$.
Let $c$ be a real number such that $||Y+c||_{r}>0$ and $\xi\in(0,\,||Y+c||_{r}]$.
Then there exists a random variable $Y^{*}$ such that: 
\begin{enumerate}
\item $P_{Y^{*}}=P_{Y}$ and $Y^{*}$ is independent of $X$, 
\item $P(|Y^{*}-Y|>\xi)\leq11\left(\xi^{-1}||Y+c||_{r}\right)^{r/(2r+1)}\times\left[\alpha\left(\sigma(X),\,\sigma(Y)\right)\right]^{2r/(2r+1)}$. 
\end{enumerate}
\end{lem}

\subsection*{Spatial block decomposition}

$ $\\

Let $Y_{\mathbf{u}}=\zeta_{\mathbf{v}=([u_{i}]+1,\,1\leq i\leq N)}$,
$\mathbf{u}\in\mathbb{R}^{N}$. The following spatial blocking idea
here is that of Tran \cite{tran} and Politis and Romano \cite{polotisromano}.

Let $\Delta_{\mathbf{i}}=\int_{(i_{1}-1)}^{i_{1}}...\int_{(i_{N}-1)}^{i_{N}}Y_{\mathbf{u}}d\mathbf{u}$
. Then,

\[
S_{\mathbf{n}}=\int_{0}^{n_{1}}...\int_{0}^{n_{N}}Y_{\mathbf{u}}d\mathbf{u}=\sum_{\begin{array}{c}
1\leq i_{k}\leq n_{k}\\
k=1,...N\end{array}}\Delta_{\mathbf{i}}.\]

So, $S_{\mathbf{n}}$ is the sum of $2^{N}P^{N}$ $q_{1}\times q_{2}\times\cdots\times q_{N}$
terms $\Delta_{\mathbf{i}}$. And each of them is an integral of $Y_{\mathbf{u}}$
over a cubic block of side $p$. Let consider the classical block
decomposition:

\[
U(1,\mathbf{n},\mathbf{j})=\sum_{{k_{i}=2j_{i}p+1},\,\,{1\leq i\leq N}}^{(2j_{i}+1)p}\Delta_{\mathbf{k}},\]

\[
U(2,\mathbf{n},\mathbf{j})=\sum_{k_{i}=2j_{i}p+1,\,\,1\leq i\leq N-1}^{(2j_{i}+1)p}\,\,\,\sum_{k_{N}=(2j_{N}+1)p+1}^{2(j_{N}+1)p}\Delta_{\mathbf{k}},\]

\[
U(3,\mathbf{n},x,\mathbf{j})=\sum_{k_{i}=2j_{i}p+1,\,\,1\leq i\leq N-2}^{(2j_{i}+1)p}\,\,\,\sum_{k_{N-1}=(2j_{N-1}+1)p+1}^{2(j_{N-1}+1)p}\,\,\,\sum_{k_{N}=2j_{N}p+1}^{(2j_{N}+1)p}\Delta_{\mathbf{k}},\]

\[
U(4,\mathbf{n},\mathbf{j})=\sum_{k_{i}=2j_{i}p+1,\,\,1\leq i\leq N-2}^{(2j_{i}+1)p}\,\,\,\sum_{k_{N-1}=(2j_{N-1}+1)p+1}^{2(j_{N-1}+1)p}\,\,\,\sum_{k_{N}=(2j_{N}+1)p+1}^{2(j_{N}+1)p}\Delta_{\mathbf{k}},\]
 and so on. Note that

\[
U(2^{N-1},\mathbf{n},\mathbf{j})=\sum_{k_{i}=(2j_{i}+1)p+1,\,\,1\leq i\leq N-1}^{2(j_{i}+1)p}\,\,\,\sum_{k_{N}=2j_{N}p+1}^{(2j_{N}+1)p}\Delta_{\mathbf{k}}.\]

Finally, \[
U(2^{N},\mathbf{n},\mathbf{j})=\sum_{k_{i}=(2j_{i}+1)p+1,\,\,1\leq i\leq N}^{2(j_{i}+1)p}\Delta_{\mathbf{k}}.\]

So, \begin{equation}
S_{\mathbf{n}}=\sum_{i=1}^{2^{N}}T(\mathbf{n},i),\label{eq:gddemo1}\end{equation}
 with $T(\mathbf{n},i)=\sum_{j_{l}=0,\, l=1,...,N}^{q_{l}-1}U(i,\mathbf{n},\mathbf{j})$.

If $n_{i}\neq2pt_{i}$, $i=1,...,N$, for all set of integers $t_{1},...,t_{N}$,
then a term, say $T\left(\mathbf{n},\,2^{N}+1\right)$ containing
all the $\Delta_{\mathbf{k}}$'s at the end, and not included in the
blocks above, can be added (see Tran \cite{tran} or Biau and Cadre
\cite{biau}). This extra term does not change the result of previous
proof.

\subsection*{Proof of Lemma \ref{lem:Dabo-Yao06}. }

$ $\\

Using (\ref{eq:gddemo1}) it suffices to show that

\begin{equation}
\mathbf{P}\left(|T(\mathbf{n},i)|>\frac{\hat{\mathbf{n}}\varepsilon}{2^{N}}\right)\leq2\,\textrm{exp}\left(-\frac{\varepsilon^{2}}{4v^{2}(\mathbf{q)}}\hat{\mathbf{q}}\right)+\hat{\mathbf{q}}\times11\left(1+\frac{4C\, p^{N}M_{2}^{2/r}}{\varepsilon}\right)^{r/(2r+1)}\alpha([p])^{2r/(2r+1)}\label{eq:Tnigd}\end{equation}
 for each $1\leq i\leq2^{N}.$

Without loss of generality we will show (\ref{eq:Tnigd}) for $i=1$.
Now, we enumerate (as it is often done in this case) in arbitrary
way the $\hat{\mathbf{q}}=q_{1}\times q_{2}\times\cdots\times q_{N}$
terms $U(1,\mathbf{n},\mathbf{j})$ of sum of $T(\mathbf{n},1)$ that
we call $W_{1},...,W_{\hat{\mathbf{q}}}$. Note that the $U(1,\mathbf{n},\mathbf{j})$
are measurable with respect to the $\sigma-$field generated by $Y_{\mathbf{u}}$
with $\mathbf{u}$ such that $2j_{i}p\leq u_{i}\leq(2j_{i}+1)p$,
$i=1,...,N$.

These sets of sites are separated by a distance at least $p$ and
since for all $m=1,...,\mathbf{\hat{\mathbf{q}}}$ there exists a
$\mathbf{j}(m)$ such that $W_{m}=U(1,\mathbf{n,j}(m))$ which have
the same distribution as $W_{m}^{*}$ ,\[
\mathbf{E}|W_{m}|^{r}=\mathbf{E}|W_{m}^{*}|^{r}=\mathbf{E}\left|\int_{2j_{1}(m)p}^{(2j_{1}(m)+1)p}...\int_{2j_{N}(m)p}^{(2j_{N}(m)+1)p}Y_{\mathbf{u}}d\mathbf{u}\right|^{r},\, r\in[1,\,+\infty].\]
 Noting that \begin{align*}
\int_{2j_{k}(m)p}^{(2j_{k}(m)+1)p}Y_{\mathbf{u}}\, d\mathbf{u} & =\int_{2j_{k}(m)p}^{[2j_{k}(m)p]+1}Y_{\mathbf{u}}\, d\mathbf{u}+\sum_{v_{k}=[2j_{k}(m)p]+2}^{[(2j_{k}(m)+1)p]}\zeta_{\mathbf{v}}+\int_{[(2j_{k}(m)+1)p]}^{2j_{k}(m)+1)p}Y_{\mathbf{u}}\, d\mathbf{u}\\
\\ & =\left([2j_{k}(m)p]+1-2j_{k}(m)p\right)\zeta_{(\mathbf{v},\, v_{k}=[2j_{k}(m)p]+1)}+\sum_{v_{k}=[2j_{k}(m)p]+2}^{[(2j_{k}(m)+1)p]}\zeta_{\mathbf{v}}\\
+ & \left((2j_{k}(m)+1)p-[(2j_{k}(m)+1)p]\right)\zeta_{(\mathbf{v},\, v_{k}=[(2j_{k}(m)+1)p]+1)}\\
 & =\sum_{v_{k}=[2j_{k}(m)p]+1}^{[(2j_{k}(m)+1)p]+1}w(\mathbf{j,v})_{k}\,\zeta_{\mathbf{v}}\end{align*}
 and $|w(\mathbf{j,v})_{k}|\leq1$ $\forall k=1,...,N$, we have by
using Minkovski's inequality and \ref{eq:hyplemdabo-yao} one get
\begin{equation}
\mathbf{E}\left|\,\frac{W_{m}}{p^{N}}\,\right|^{r}\leq c^{r-2}r!\, M_{2}^{2}\,,\forall r\geq2.\label{eq:momentWnborn}\end{equation}
 Then, using recursively the version of Bradley's lemma gives in Lemma~\ref{lem:(-Bosq,-1997)}
we define independent random variables $W_{1}^{*},...,W_{\hat{\mathbf{q}}}^{*}$
such that for all $r\in[1,+\infty]$ and for all $m=1,...,\hat{\mathbf{q}}$,
$W_{m}^{*}$ has the same distribution with $W_{m}$ and setting $\omega_{r}^{r}=p^{rN}c^{r-2}M_{2}^{2}$,
we have:\[
P(|W_{m}-W_{m}^{*}|>\xi)\leq11\left(\frac{||W_{m}+\omega_{r}||_{r}}{\xi}\right)^{r/(2r+1)}\alpha([p])^{2r/(2r+1)},\]
 where, $c=\delta\omega_{r}p$ and $\xi=\min\left(\frac{\hat{\mathbf{n}}\varepsilon}{2^{N+1}\hat{\mathbf{q}}},\,(\delta-1)\omega_{r}p^{N}\right)=\min\left(\frac{\varepsilon p^{N}}{2},\,(\delta-1)\omega_{r}p^{N}\right)$
for some $\delta>1$ specified below. Note that for each $m$, \[
||W_{m}+c||_{r}\geq c-||W_{m}||_{r}\geq(\delta-1)\omega_{r}p^{N}>0\]
 so that $0<\xi<||W_{m}+c||_{r}$ as required in Lemma \ref{lem:(-Bosq,-1997)}.

Then, if $\delta=1+\frac{\varepsilon}{2\omega_{r}}$,\[
P(|W_{m}-W_{m}^{*}|>\xi)\leq11\left(1+\frac{4\omega_{r}}{\varepsilon}\right)^{r/(2r+1)}\alpha([p])^{2r/(2r+1)}\]

and \[
P\left(\sum_{m=1}^{\hat{\mathbf{q}}}\left|W_{m}-W_{m}^{*}\right|>\frac{\hat{\mathbf{n}}\varepsilon}{2^{N+1}}\right)\leq\hat{\mathbf{q}}\times11\left(1+\frac{4\omega_{r}}{\varepsilon}\right)^{r/(2r+1)}\alpha([p])^{2r/(2r+1)}.\]

Now, note that Inequality (\ref{eq:momentWnborn}) also leads (by
Bernstein's inequality) to :

\[
P\left(\left|\sum_{m=1}^{\hat{\mathbf{q}}}W_{m}^{*}\right|>\frac{\hat{\mathbf{n}}\varepsilon}{2^{N+1}}\right)\leq2\,\textrm{exp}\left(-\frac{\left(\frac{\hat{\mathbf{n}}\varepsilon}{2^{N+1}}\right)^{2}}{4\sum_{m=1}^{\hat{\mathbf{q}}}\mathbf{E}\, W_{m}^{2}+\frac{c\hat{\mathbf{n}}p^{N}}{2^{N+1}}\varepsilon}\right)\]
 Thus \[
\begin{array}{ccc}
P(|T(\mathbf{n},1)|>\frac{\hat{\mathbf{n}}\varepsilon}{2^{N}}) & \leq & 2\textrm{exp}\left(-\frac{\hat{\mathbf{q}}\varepsilon^{2}}{4(\, M_{2}^{2}+2^{N}c\varepsilon)}\right)+\hat{\mathbf{q}}\times11\left(1+\frac{4c\, p^{N}M_{2}^{2/r}}{\varepsilon}\right)^{r/(2r+1)}\alpha([p])^{2r/(2r+1)}\end{array}\]

Then, since $\hat{\mathbf{q}}=q_{1}\times...\times q_{N}$ and $\hat{\mathbf{n}}=2^{N}p^{N}\hat{\mathbf{q}}$,
we get inequality (\ref{eq:Tnigd}) the proof is completed by noting
that $P(|S_{\mathbf{n}}|>\hat{\mathbf{n}}\varepsilon)\leq2^{N}P(|T(\mathbf{n},i)|>\frac{\hat{\mathbf{n}}\varepsilon}{2^{N}})$.~~~~~~~~~~$\square$

\subsection{Proof of the Theorem \ref{thm:conv-prob} }

We will prove the desired result on $\Sigma_{e,\,\mathbf{n}}-\Sigma_{e}$
using an intermediate matrix \[
\overline{\Sigma}_{e,\mathbf{n}}=\frac{1}{\hat{\mathbf{n}}}\sum_{\mathbf{i}\in\mathcal{I}_{\mathbf{n}}}r(Y_{\mathbf{i}})r(Y_{\mathbf{i}})^{T}.\]
 Start with the following decomposition \[
\Sigma_{e,\,\mathbf{n}}-\Sigma_{e}=\Sigma_{e,\,\mathbf{n}}-\overline{\Sigma}_{e,\mathbf{n}}+\overline{\Sigma}_{e,\mathbf{n}}-\Sigma_{e}.\]
 We first show that: \begin{equation}
\Sigma_{e,\,\mathbf{n}}-\overline{\Sigma}_{e,\,\mathbf{n}}=\mathcal{O}_{p}\left(\frac{1}{\hat{\mathbf{n}}^{\frac{1}{2}+\delta}}+\frac{\Psi_{\mathbf{n}}^{2}}{e_{\mathbf{n}}^{2}}\right).\label{eq:Sigmae1}\end{equation}
 To this aim, we set~: \begin{equation}
\Sigma_{e,\,\mathbf{n}}-\overline{\Sigma}_{e,\,\mathbf{n}}=S_{\mathbf{n},\,1}+S_{\mathbf{n},\,2}+S_{\mathbf{n},\,3}\label{eq:Sigmae2}\end{equation}

with \[
S_{\mathbf{n},1}=\frac{1}{\hat{\mathbf{n}}}\sum_{\mathbf{i}\in\mathcal{I}_{\mathbf{n}}}\left(\hat{r}_{e_{\mathbf{n}}}(Y_{\mathbf{i}})-r(Y_{\mathbf{i}})\right)\left(\hat{r}_{e_{\mathbf{n}}}(Y_{\mathbf{i}})-r(Y_{\mathbf{i}})\right)^{T},\]
 \[
S_{\mathbf{n},\,2}=\frac{1}{\hat{\mathbf{n}}}\sum_{\mathbf{i}\in\mathcal{I}_{\mathbf{n}}}r(Y_{\mathbf{i}})\,\left(\hat{r}_{e_{n}}(Y_{i})-r(Y_{\mathbf{i}})\right)^{T}\]
 and

\[
S_{\mathbf{n},\,3}=\frac{1}{\hat{\mathbf{n}}}\sum_{\mathbf{i}\in\mathcal{I}_{\mathbf{n}}}\left(\hat{r}_{e_{n}}(Y_{i})-r(Y_{\mathbf{i}})\right)\, r(Y_{\mathbf{i}})^{T}.\]
 Note that $S_{\mathbf{n},\,3}^{T}=S_{\mathbf{n},\,2}$, hence we
only need to control the rate of convergence of the first two terms
$S_{\mathbf{n},\,1}$ and $S_{\mathbf{n},\,2}$ \vskip .1in We will
successively prove that \[
S_{\mathbf{n},1}=\mathcal{O}_{p}\left(\frac{\Psi_{\mathbf{n}}^{2}}{e_{\mathbf{n}}^{2}}\right),\]
 and \[
S_{\mathbf{n},2}=\mathcal{O}_{p}\left(\frac{\Psi_{\mathbf{n}}^{2}}{e_{\mathbf{n}}^{2}}+h_{\mathbf{n}}^{k}\right)\]
 this latter will immediately implies that \[
S_{\mathbf{n},3}=\mathcal{O}_{p}\left(\frac{\Psi_{\mathbf{n}}^{2}}{e_{\mathbf{n}}^{2}}+h_{\mathbf{n}}^{k}\right).\]

\begin{itemize}
\item Control on $S_{\mathbf{n},\,1}$ 
\end{itemize}
Since for each $y\in\mathbb{R}$ : \begin{equation}
\hat{r}_{e_{\mathbf{n}}}(y)-r(y)=\frac{r(y)}{f_{e_{\mathbf{n}}}(y)}\left(f(y)-f_{e_{\mathbf{n}}}(y)\right)+\frac{1}{\hat{f}_{e_{\mathbf{n}}}(y)}\left(\varphi_{\mathbf{n}}(y)-\varphi(y)\right)\label{eq:rhat-r}\end{equation}
 and\begin{equation}
f(y)-f_{e_{\mathbf{n}}}(y)=f(y)-f_{\mathbf{n}}(y)+(f_{\mathbf{n}}(y)-e_{\mathbf{n}})\mathbf{1}_{\{f_{\mathbf{n}}(y)<e_{\mathbf{n}}\}},\label{eq:f-fen}\end{equation}
 for each $\mathbf{i}\in(\mathbb{N}^{*})^{N}$ \[
\left\Vert r_{e_{n}}(Y_{\mathbf{i}})-r(Y_{\mathbf{i}})\right\Vert \leq\frac{\left\Vert r(Y_{\mathbf{i}})\right\Vert }{e_{\mathbf{n}}}\,||f_{\mathbf{n}}-f||_{\infty}+2\left\Vert r(Y_{\mathbf{i}})\right\Vert \mathbf{1}_{\{f_{\mathbf{n}}(Y_{\mathbf{i}})<e_{\mathbf{n}}\}}+\frac{\left\Vert \varphi_{\mathbf{n}}-\varphi\right\Vert _{\infty}}{e_{\mathbf{n}}}.\]

and \[
\left\Vert r_{e_{\mathbf{n}}}(Y_{\mathbf{i}})-r(Y_{\mathbf{i}})\right\Vert ^{2}\leq3\left[\left\Vert r(Y_{\mathbf{i}})\right\Vert ^{2}\frac{||f_{\mathbf{n}}-f||_{\infty}^{2}}{e_{n}^{2}}\,+4\left\Vert r(Y_{\mathbf{i}})\right\Vert ^{2}\mathbf{1}_{\{f_{\mathbf{n}}(Y_{\mathbf{i}})<e_{\mathbf{n}}\}}+\frac{||\varphi_{\mathbf{n}}-\varphi||_{\infty}^{2}}{e_{\mathbf{n}}^{2}}\right].\]

Using the following inequality (see Ferré and Yao \cite{Ferreyao}
for details): \begin{equation}
\mathbf{1}_{\{f_{\mathbf{n}}(Y_{\mathbf{i}})<e_{n}\}}\leq\mathbf{1}_{\{f(Y_{\mathbf{i}})<e_{\mathbf{n}}\}}+\frac{||f_{\mathbf{n}}-f||_{\infty}^{2}}{e_{n}^{2}},\label{eq:ferre-yao}\end{equation}
 and by results on Lemmas \ref{lem: gLaw-of-large} and \ref{thm:Carbon-et-alPro},
we have: \[
S_{\mathbf{n},1}\leq\frac{C}{\hat{\mathbf{n}}}\sum_{\mathbf{i}\in\mathcal{I}_{\mathbf{n}}}\left\Vert r(Y_{\mathbf{i}})\right\Vert ^{2}\mathbf{1}_{\{f(Y_{\mathbf{i}})<e_{\mathbf{n}}\}}+\mathcal{O}_{p}\left(\frac{\Psi_{\mathbf{n}}^{2}}{e_{\mathbf{n}}^{2}}\right),\, C>0.\]
 Now, noting that \[
\frac{1}{\hat{\mathbf{n}}}\sum_{\mathbf{i}\in\mathcal{I}_{\mathbf{n}}}\left\Vert r(Y_{\mathbf{i}})\right\Vert ^{2}\mathbf{1}_{\{f(Y_{\mathbf{i}})<e_{\mathbf{n}}\}}\leq e_{\mathbf{n}}^{2}\,\frac{1}{\hat{\mathbf{n}}}\sum_{\mathbf{i}\in\mathcal{I}_{\mathbf{n}}}\frac{\left\Vert r(Y_{\mathbf{i}})\right\Vert ^{2}}{f(Y_{\mathbf{i}})^{2}}\mathbf{1}_{\{f(Y_{\mathbf{i}})<e_{\mathbf{n}}\}},\]
 we have (since $\mathbf{E}\left(\frac{||r(Y_{\mathbf{i}})||^{2}}{f(Y_{\mathbf{i}})^{2}}\mathbf{1}_{\{f(Y_{\mathbf{i}})<e_{\mathbf{n}}\}}\right)=\mathcal{O}\left(\frac{1}{\hat{\mathbf{n}}^{1+\delta}}\right)$
by assumption): \begin{equation}
\frac{1}{\hat{\mathbf{n}}}\sum_{\mathbf{i}\in\mathcal{I}_{\mathbf{n}}}\left\Vert r(Y_{\mathbf{i}})\right\Vert ^{2}\mathbf{1}_{\{f(Y_{\mathbf{i}})<e_{\mathbf{n}}\}}=\mathcal{O}_{p}\left(\frac{e_{\mathbf{n}}^{2}}{\hat{\mathbf{n}}^{1+\delta}}\right)\label{eq:meansqen}\end{equation}
 and \[
S_{\mathbf{n},1}=\mathcal{O}_{p}\left(\frac{e_{\mathbf{n}}^{2}}{\hat{\mathbf{n}}^{1+\delta}}+\frac{\Psi_{\mathbf{n}}^{2}}{e_{\mathbf{n}}^{2}}\right)\]
 because of Assumption $\mathbf{E}\left(\frac{\left\Vert r(Y)\right\Vert ^{2}}{f(Y)^{2}}\mathbf{1}_{\{f(Y)<e_{\mathbf{n}}\}}\right)=\mathcal{O}\left(\frac{1}{\hat{\mathbf{n}}^{1+\delta}}\right)$.

Now, since $\Psi_{\mathbf{n}}=h_{\mathbf{n}}^{k}+\sqrt{\frac{\log\hat{\mathbf{n}}}{\hat{\mathbf{n}}\, h_{\mathbf{n}}}}$
and $\frac{e_{\mathbf{n}}}{\hat{\mathbf{n}}^{\frac{1+\delta}{2}}}\leq C\,\sqrt{\frac{\log\hat{\mathbf{n}}}{\hat{\mathbf{n}}\, h_{\mathbf{n}}}}$
(for $\hat{\mathbf{n}}$ large and $C>0$ an arbitrary constante),
we have: \begin{equation}
S_{\mathbf{n},1}=\mathcal{O}_{p}\left(\frac{\Psi_{\mathbf{n}}^{2}}{e_{\mathbf{n}}^{2}}\right).\label{eq:Sn1}\end{equation}

\begin{itemize}
\item Control on $S_{\mathbf{n},\,2}$ . 
\end{itemize}
Noting that : $\frac{1}{f_{e_{n}}}=\frac{1}{f}+\frac{f-\tilde{f}_{e_{\mathbf{n}}}}{\tilde{f}_{e_{\mathbf{n}}}f}+\frac{\tilde{f}_{e_{\mathbf{n}}}-f_{e_{\mathbf{n}}}}{\tilde{f}_{e_{\mathbf{n}}}\hat{f}_{e_{\mathbf{n}}}}$$=\frac{1}{f}+\frac{f-e_{\mathbf{n}}}{\tilde{f}_{e_{\mathbf{n}}}f}\mathbf{1}_{\{f<e_{\mathbf{n}}\}}+\frac{\tilde{f}_{e_{\mathbf{n}}}-f_{e_{\mathbf{n}}}}{\tilde{f}_{e_{\mathbf{n}}}f_{e_{\mathbf{n}}}}$,
with $\tilde{f}_{e_{\mathbf{n}}}=\max\{f,e_{\mathbf{n}}\}$, we have:

\begin{eqnarray*}
S_{\mathbf{n},\,2} & = & \frac{1}{\hat{\mathbf{n}}}\sum_{\mathbf{i}\in\mathcal{I}_{\mathbf{n}}}\frac{r(Y_{\mathbf{i}})r(Y_{\mathbf{i}})^{T}}{f_{e_{\mathbf{n}}}(Y_{\mathbf{i}})}\left(f(Y_{\mathbf{i}})-f_{e_{\mathbf{n}}}(Y_{\mathbf{i}})\right)+\frac{r(Y_{\mathbf{i}})}{f_{e_{\mathbf{n}}}(Y_{\mathbf{i}})}\left(\varphi_{\mathbf{n}}(Y_{\mathbf{i}})-\varphi(Y_{\mathbf{i}})\right)^{T}\\
 & = & \frac{1}{\hat{\mathbf{n}}}\sum_{\mathbf{i}\in\mathcal{I}_{\mathbf{n}}}\frac{r(Y_{\mathbf{i}})r(Y_{\mathbf{i}})^{T}}{f(Y_{\mathbf{i}})}\left(f(Y_{\mathbf{i}})-f_{\mathbf{n}}(Y_{\mathbf{i}})\right)+\frac{r(Y_{\mathbf{i}})}{f(Y_{\mathbf{i}})}\left(\varphi_{\mathbf{n}}(Y_{\mathbf{i}})-\varphi(Y_{\mathbf{i}})\right)^{T}+R_{\mathbf{n}_{1}}+R_{\mathbf{n}_{2}}.\end{eqnarray*}
 where\begin{eqnarray*}
R_{\mathbf{n}_{1}}(Y_{\mathbf{i}}) & = & r(Y_{\mathbf{i}})\,\left[r(Y_{\mathbf{i}})^{T}\left(f(Y_{\mathbf{i}})-f_{\mathbf{n}}(Y_{\mathbf{i}})\right)+\left(\varphi_{\mathbf{n}}(Y_{\mathbf{i}})-\varphi(Y_{\mathbf{i}})\right)^{T}\right]\\
 &  & \left(\frac{1}{f(Y_{\mathbf{i}})}\mathbf{1}_{\{f(Y_{\mathbf{i}})<e_{\mathbf{n}}\}}+\frac{\tilde{f}_{e_{\mathbf{n}}}(Y_{\mathbf{i}})-f_{e_{\mathbf{n}}}(Y_{\mathbf{i}})}{\tilde{f}_{e_{\mathbf{n}}}(Y_{\mathbf{i}})f_{e_{\mathbf{n}}}(Y_{\mathbf{i}})}\right)\end{eqnarray*}
 and \[
R_{\mathbf{n}_{2}}=\frac{r(Y_{\mathbf{i}})r(Y_{\mathbf{i}})^{T}}{f_{e_{\mathbf{n}}}(Y_{\mathbf{i}})}\left(f_{\mathbf{n}}(Y_{i})-f_{e_{\mathbf{n}}}(Y_{\mathbf{i}})\right).\]
Futhermore :

\begin{itemize}
\item since for all $y\in\mathbb{R}$ we have $\frac{1}{\tilde{f}_{e_{\mathbf{n}}}(y)f_{e_{\mathbf{n}}}(y)}\leq\frac{1}{e_{\mathbf{n}}^{2}}$
and by several calculus we also have $\left|\tilde{f}_{e_{\mathbf{n}}}(y)-f_{e_{\mathbf{n}}}(y)\right|\leq\left|f(y)-f_{\mathbf{n}}(y)\right|$
and then $||\tilde{f}_{e_{\mathbf{n}}}-f_{e_{\mathbf{n}}}||_{\infty}\leq||f_{\mathbf{n}}-f||_{\infty}$
, we also have one hand: {\footnotesize \begin{eqnarray}
\,\,\,\,\,\:\, R_{\mathbf{n}_{1}} & \leq & \frac{1}{\hat{\mathbf{n}}}\sum_{\mathbf{i}\in\mathcal{I}_{\mathbf{n}}}\left(||r(Y_{\mathbf{i}})||\,||\varphi_{\mathbf{n}}-\varphi||_{\infty}+||r(Y_{\mathbf{i}})||^{2}\,||f_{\mathbf{n}}-f||_{\infty}\right)\,\left(\frac{1}{f(Y_{\mathbf{i}})}\mathbf{1}_{\{f<e_{\mathbf{n}}\}}+\frac{||f_{\mathbf{n}}-f||_{\infty}}{e_{\mathbf{n}}^{2}}\right)\label{eq:Rn1}\end{eqnarray}
}{\footnotesize \par}
\item on the other hand we have \begin{eqnarray*}
R_{\mathbf{n}_{2}} & \leq & \frac{1}{\hat{\mathbf{n}}}\sum_{\mathbf{i}\in\mathcal{I}_{\mathbf{n}}}||r(Y_{\mathbf{i}})||^{2}\,\frac{|f_{\mathbf{n}}(Y_{\mathbf{i}})-e_{\mathbf{n}}|}{f_{e_{\mathbf{n}}}(Y_{\mathbf{i}})}\mathbf{1}_{\{f_{\mathbf{n}}(Y_{i})<e_{\mathbf{n}}\}}\\
 & \leq & \frac{2}{\hat{\mathbf{n}}}\sum_{\mathbf{i}\in\mathcal{I}_{\mathbf{n}}}||r(Y_{\mathbf{i}})||^{2}\,\mathbf{1}_{\{f_{\mathbf{n}}(Y_{i})<e_{\mathbf{n}}\}}.\end{eqnarray*}
because for all $y\in\mathbb{R}$ , $|f_{\mathbf{n}}(y)-f_{e_{\mathbf{n}}}(y)|=\left|f_{\mathbf{n}}(y)-e_{\mathbf{n}}\right|\mathbf{1}_{\{f_{\mathbf{n}}(y)<e_{\mathbf{n}}\}}\leq2e_{\mathbf{n}}\mathbf{1}_{\{f_{\mathbf{n}}(y)<e_{\mathbf{n}}\}}$.
\end{itemize}
Then, it follows from (\ref{eq:ferre-yao} and \ref{eq:meansqen})
that:

\[
R_{\mathbf{n}_{2}}=\mathcal{O}_{p}\left(\frac{e_{\mathbf{n}}^{2}}{\hat{\mathbf{n}}^{1+\delta}}+\frac{\Psi_{\mathbf{n}}^{2}}{e_{\mathbf{n}}^{2}}\right)\]
as for $S_{\mathbf{n}_{1}}$, we deduce: \[
R_{\mathbf{n}_{2}}=\mathcal{O}_{p}\left(\frac{\Psi_{\mathbf{n}}^{2}}{e_{\mathbf{n}}^{2}}\right)\]

Now, observious that,

\[
\frac{1}{\hat{\mathbf{n}}}\sum_{\mathbf{i}\in\mathcal{I}_{\mathbf{n}}}\frac{||r(Y_{\mathbf{i}})||^{2}}{f(Y_{\mathbf{i}})}\mathbf{1}_{\{f(Y_{\mathbf{i}})<e_{\mathbf{n}}\}}\leq e_{\mathbf{n}}\,\frac{1}{\hat{\mathbf{n}}}\sum_{\mathbf{i}\in\mathcal{I}_{\mathbf{n}}}\frac{||r(Y_{\mathbf{i}})||^{2}}{f(Y_{\mathbf{i}})^{2}}\mathbf{1}_{\{f(Y_{\mathbf{i}})<e_{\mathbf{n}}\}},\]
 we have (as previously): \begin{equation}
\frac{1}{\hat{\mathbf{n}}}\sum_{\mathbf{i}\in\mathcal{I}_{\mathbf{n}}}\frac{||r(Y_{\mathbf{i}})||^{2}}{f(Y_{\mathbf{i}})}\mathbf{1}_{\{f(Y_{\mathbf{i}})<e_{\mathbf{n}}\}}=\mathcal{O}_{p}\left(\frac{e_{\mathbf{n}}}{\hat{\mathbf{n}}^{1+\delta}}\right).\label{eq:Rn11}\end{equation}
 Moreover, since $\mathbf{E}\left(\frac{\left\Vert r(Y)\right\Vert ^{2}}{f(Y)^{2}}\mathbf{1}_{\{f(Y)<e_{\mathbf{n}}\}}\right)=\mathcal{O}\left(\frac{1}{\hat{\mathbf{n}}^{1+\delta}}\right)$,
we also have: \begin{equation}
\frac{1}{\hat{\mathbf{n}}}\sum_{\mathbf{i}\in\mathcal{I}_{\mathbf{n}}}\frac{||r(Y_{\mathbf{i}})||}{f(Y_{\mathbf{i}})}\mathbf{1}_{\{f(Y_{\mathbf{i}})<e_{\mathbf{n}}\}}=\mathcal{O}_{p}\left(\frac{1}{\hat{\mathbf{n}}^{\frac{1+\delta}{2}}}\right)\label{eq:Rn12}\end{equation}

So combining (\ref{eq:Rn1}), (\ref{eq:Rn11}) and (\ref{eq:Rn12}),
we get: \[
R_{\mathbf{n}_{1}}=\mathcal{O}_{p}\left(\frac{e_{\mathbf{n}}\Psi_{\mathbf{n}}}{\hat{\mathbf{n}}^{1+\delta}}+\frac{\Psi_{\mathbf{n}}}{\hat{\mathbf{n}}^{\frac{1+\delta}{2}}}+\frac{\Psi_{\mathbf{n}}^{2}}{e_{\mathbf{n}}^{2}}\right)=\mathcal{O}_{p}\left(\frac{\Psi_{\mathbf{n}}}{\hat{\mathbf{n}}^{\frac{1+\delta}{2}}}+\frac{\Psi_{\mathbf{n}}^{2}}{e_{\mathbf{n}}^{2}}\right);\]
 and since $\frac{e_{\mathbf{n}}}{\hat{\mathbf{n}}^{\frac{1+\delta}{2}}}\leq C\,\sqrt{\frac{\log\hat{\mathbf{n}}}{\hat{\mathbf{n}}\, h_{\mathbf{n}}}}$
(for $\hat{\mathbf{n}}$ large) we have: \[
R_{\mathbf{n}_{1}}=\mathcal{O}_{p}\left(\frac{\Psi_{\mathbf{n}}^{2}}{e_{\mathbf{n}}^{2}}\right).\]
Then,

\[
S_{\mathbf{n},\,2}=S_{\mathbf{n},\,2}^{(1)}+S_{\mathbf{n},\,2}^{(2)}+\mathcal{O}_{p}\left(\frac{\Psi_{\mathbf{n}}^{2}}{e_{\mathbf{n}}^{2}}\right);\]

with \[
S_{\mathbf{n},\,2}^{(1)}=\frac{1}{\hat{\mathbf{n}}}\sum_{i=1}\frac{r(Y_{\mathbf{i}})r(Y_{\mathbf{i}})^{T}}{f(Y_{\mathbf{i}})}\left(f_{\mathbf{n}}(Y_{\mathbf{i}})-f(Y_{\mathbf{i}})\right)\]
 and \[
S_{\mathbf{n},\,2}^{(2)}=\frac{1}{\hat{\mathbf{n}}}\sum_{i=1}\frac{r(Y_{\mathbf{i}})}{f(Y_{\mathbf{i}})}\left(\varphi_{\mathbf{n}}(Y_{\mathbf{i}})-\varphi(Y_{\mathbf{i}})\right)^{T}.\]
 To finish, we are going to show that \[
S_{\mathbf{n},\,2}^{(1)}=\mathcal{O}_{p}(h_{\mathbf{n}}^{k}+\frac{1}{\hat{\mathbf{n}}h_{\mathbf{n}}})\]
 \[
S_{\mathbf{n},\,2}^{(2)}=\mathcal{O}_{p}(h_{\mathbf{n}}^{k}+\frac{1}{\hat{\mathbf{n}}h_{\mathbf{n}}})\]

Note that: \[
S_{\mathbf{n},\,2}^{(1)}=\frac{1}{\hat{\mathbf{n}}}\sum_{i=1}\tau(Y_{\mathbf{i}})\, f(Y_{\mathbf{i}})-\frac{1}{h_{\mathbf{n}}}V_{\mathbf{n}}\]
 where $\tau(.)$ is a function defined by $\tau(y)=\frac{r(y)\, r(y)^{T}}{f(y)\,}$
for $y\in\mathbb{R}$ and \[
V_{\mathbf{n}}=\frac{1}{\hat{\mathbf{n}}^{2}}\sum_{\mathbf{i},\mathbf{j}\in\mathcal{I}_{\mathbf{n}}}\tau(Y_{\mathbf{i}})\, K_{h_{\mathbf{n}}}(Y_{\mathbf{i}}-Y_{\mathbf{j}})\]
 is a second-order Von Mises functional statistic which associated
U-statistic is: \[
U_{\mathbf{n}}=\frac{1}{2\hat{\mathbf{n}}(\hat{\mathbf{n}}-1)}\sum_{i=1}\sum_{j\neq i}[\tau(Y_{\mathbf{i}})+\tau(Y_{\mathbf{j}})]K_{h_{\mathbf{n}}}(Y_{\mathbf{i}}-Y_{\mathbf{j}}).\]

Since: $V_{\mathbf{n}}=U_{\mathbf{n}}+\mathcal{O}_{p}(\frac{1}{\hat{\mathbf{n}}})$,
\[
S_{\mathbf{n},\,2}^{(1)}=\frac{1}{\hat{\mathbf{n}}}\sum_{i=1}\tau(Y_{\mathbf{i}})\, f(Y_{\mathbf{i}})-\frac{1}{h_{\mathbf{n}}}U_{\mathbf{n}}+\mathcal{O}_{p}\left(\frac{1}{\hat{\mathbf{n}}h_{\mathbf{n}}}\right).\]
 We apply Lemma \ref{lem:Ustat-pro} with, $m=2$, $h(y_{1},y_{2})=\left[\tau(y_{1})+\tau(y_{2})\right]K_{h_{\mathbf{n}}}(y_{1}-y_{2})$
\[
h_{1}(y)=\frac{1}{2}\left[\tau(y).f*K_{h_{\mathbf{n}}}(y)+(\tau.f)*K_{h_{n}}(y)\right],\]
 and

\[
\Theta(F)=\mathbf{E}\left(h_{1}(Y)\right)=\mathbf{E}\,\left(\tau(y).f*K_{h_{\mathbf{n}}}(y)\right).\]
 Since \[
||h(Y_{1},Y_{2})||_{4+\delta}\leq C.||\tau(Y)||_{4+\delta}<\infty,\]
 by assumption (\ref{eq:Assump-rrT/f}) then, \[
U_{\mathbf{n}}=\Theta(F)+\frac{2}{\hat{\mathbf{n}}}\sum_{\mathbf{i}}\,\left(h_{1}(Y_{\mathbf{i}})-\Theta(F)\right)+\mathcal{O}_{p}(\frac{1}{\hat{\mathbf{n}}}).\]
 and \begin{eqnarray*}
S_{\mathbf{n},\,2}^{(1)} & = & \frac{1}{\hat{\mathbf{n}}}\sum_{i=1}\tau(Y_{\mathbf{i}})\, f(Y_{\mathbf{i}})-\frac{\Theta(F)}{h_{\mathbf{n}}}-\frac{2}{\hat{\mathbf{n}}}\sum_{\mathbf{i}}\,\left(\frac{h_{1}(Y_{\mathbf{i}})}{h_{\mathbf{n}}}-\frac{\Theta(F)}{h_{\mathbf{n}}}\right)+\mathcal{O}_{p}\left(\frac{1}{\hat{\mathbf{n}}h_{\mathbf{n}}}\right)\\
 & = & \frac{1}{\hat{\mathbf{n}}}\sum_{i=1}\tau(Y_{\mathbf{i}})\,\left(f(Y_{\mathbf{i}})-\frac{f*K_{h_{\mathbf{n}}}(Y_{\mathbf{i}})}{h_{\mathbf{n}}}\right)+\frac{\Theta(F)-(\tau.f)*K_{h_{\mathbf{n}}}(y_{\mathbf{i}})}{h_{\mathbf{n}}}+\mathcal{O}_{p}\left(\frac{1}{\hat{\mathbf{n}}h_{\mathbf{n}}}\right).\end{eqnarray*}

Since $f$ and $r(.)$ belongs to $C^{k}(\mathbb{R})$, we get, \[
\left\Vert \frac{f*K_{h_{\mathbf{n}}}(y)}{h_{\mathbf{n}}}-f(y)\right\Vert _{\infty}=\mathcal{O}(h_{\mathbf{n}}^{k})\]
 and \[
\left\Vert \frac{(\tau.f)*K_{h_{\mathbf{n}}}(y)}{h_{\mathbf{n}}}-(\tau.f)(y)\right\Vert _{\infty}=\mathcal{O}(h_{\mathbf{n}}^{k}).\]
 Then, we have \[
\frac{1}{\hat{\mathbf{n}}}\sum_{i\in\mathcal{I}_{\mathbf{n}}}\frac{\tau(Y_{\mathbf{i}})\, f*K_{h_{\mathbf{n}}}(Y_{\mathbf{i}})}{h_{\mathbf{n}}}=\frac{1}{\hat{\mathbf{n}}}\sum_{i\in\mathcal{I}_{\mathbf{n}}}(\tau.f)(Y_{\mathbf{i}})+\mathcal{O}_{p}(h_{\mathbf{n}}^{k}),\]
 \[
\frac{1}{\hat{\mathbf{n}}}\sum_{i\in\mathcal{I}_{\mathbf{n}}}\frac{(\tau.f)*K_{h_{\mathbf{n}}}(Y_{\mathbf{i}})}{h_{\mathbf{n}}}=\frac{1}{\hat{\mathbf{n}}}\sum_{i\in\mathcal{I}_{\mathbf{n}}}(\tau.f)(Y_{\mathbf{i}})+\mathcal{O}(h_{\mathbf{n}}^{k})\]
 \[
\frac{\Theta(F)}{h_{\mathbf{n}}}=\mathbf{E}((\tau.f)(Y))+\mathcal{O}(h_{\mathbf{n}}^{k}).\]
 Finally: \[
S_{\mathbf{n},\,2}^{(1)}=\mathcal{O}_{p}(h_{\mathbf{n}}^{k}+\frac{1}{\hat{\mathbf{n}}}+\frac{1}{\hat{\mathbf{n}}h_{\mathbf{n}}})=\mathcal{O}_{p}(h_{\mathbf{n}}^{k}+\frac{1}{\hat{\mathbf{n}}h_{\mathbf{n}}}).\]
 By using similar arguments and applying Lemma \ref{lem:Ustat-pro}
with $m=3$, one also gets \[
S_{\mathbf{n},\,2}^{(2)}=\mathcal{O}_{p}(h_{\mathbf{n}}^{k}+\frac{1}{\hat{\mathbf{n}}h_{\mathbf{n}}}).\]
 So, \[
S_{\mathbf{n},2}=\mathcal{O}_{p}\left(\frac{\Psi_{\mathbf{n}}^{2}}{e_{\mathbf{n}}^{2}}+h_{\mathbf{n}}^{k}+\frac{1}{\hat{\mathbf{n}}h_{\mathbf{n}}}\right)\]
 Then, equality, (\ref{eq:Sn1}), and (\ref{eq:Sn2}) lead to (\ref{eq:Sigmae1}).

Recall that $\Psi_{\mathbf{n}}=h_{\mathbf{n}}^{k}+\sqrt{\frac{\log\hat{\mathbf{n}}}{\hat{\mathbf{n}}\, h_{\mathbf{n}}}}$
. Then, the fact that there exist a real $A>0$ such that $\forall\,\hat{\mathbf{n}}>A$,
$\frac{1}{\hat{\mathbf{n}}\, h_{\mathbf{n}}}<\frac{\log\hat{\mathbf{n}}}{\hat{\mathbf{n}}\, h_{\mathbf{n}}e_{\mathbf{n}}^{2}}$
and : \begin{equation}
S_{\mathbf{n},2}=\mathcal{O}_{p}\left(\frac{\Psi_{\mathbf{n}}^{2}}{e_{\mathbf{n}}^{2}}+h_{\mathbf{n}}^{k}\right)\label{eq:Sn2}\end{equation}

Finally, using equality (\ref{eq:Sigmae1}) one has; \begin{equation}
\Sigma_{e,\,\mathbf{n}}-\Sigma_{e}=\overline{\Sigma}_{e,\mathbf{n}}-\Sigma_{e}+\mathcal{O}_{p}\left(\frac{\Psi_{\mathbf{n}}^{2}}{e_{\mathbf{n}}^{2}}+h_{\mathbf{n}}^{k}\right).\label{eq:Rel-entre-Sigma-final}\end{equation}

To complete the proof, we will use Lemma \ref{lem: gLaw-of-large}.
To this aim, it suffices to choose $\theta=\delta$ with $\delta>2N$
then $\mathbf{E}||X||^{4+\delta}<\infty$ and $\sum_{k}\alpha(k)^{\frac{\delta}{\delta+4}}<\infty$;
hence we have: \[
\overline{\Sigma}_{e,\mathbf{n}}-\Sigma_{e}=\mathcal{O}_{p}(\frac{1}{\hat{\mathbf{n}}}).\]
 which ends the proof. ~~~~~~~~~~~~~~~~~~~~~~~~~~~~~~~~~~~~~$\square$

\subsection{Proof of corollary \ref{cor:cov-prob-vitesse opt} }

$ $\\

The proof is achieved by replacing $h_{\mathbf{n}}\simeq\hat{\mathbf{n}}^{-c_{1}}$
and $e_{\mathbf{n}}\simeq\hat{\mathbf{n}}^{-c_{2}}$ with $\frac{c_{2}}{k}+\frac{1}{4k}<c_{1}<\frac{1}{2}-2c_{2}$
on equality (\ref{eq:Rel-entre-Sigma-final}) $\square$

\subsection{Proof of corrollary \ref{thm:(Central-limit-theorem)} }

$ $\\

Chosing $h_{\mathbf{n}}\simeq\hat{\mathbf{n}}^{-c_{1}}$ and $e_{\mathbf{n}}\simeq\hat{\mathbf{n}}^{-c_{2}}$
where $\frac{c_{2}}{k}+\frac{1}{4k}<c_{1}<\frac{1}{2}-2c_{2}$ on
equality (\ref{eq:Rel-entre-Sigma-final}), one gets $\Sigma_{e,\,\mathbf{n}}-\Sigma_{e}=\overline{\Sigma}_{e,\,\mathbf{n}}-\Sigma_{e}+o_{p}(\frac{1}{\sqrt{\hat{\mathbf{n}}}})$
and the central limit theorem for spatial data and Slusky's theorem
completes the proof.

\subsection*{Proof of Theorem \ref{thm:Conv-a.s} }

Let $v_{\mathbf{n}}=\left(\frac{\log\log\,\hat{\mathbf{n}}}{\hat{\mathbf{n}}}\right)^{\frac{1}{2}}$
note that since $Y$ take place on a compact set, $\frac{1}{f}$ is
bounded and replace the assumption $\left\Vert \frac{r(Y)\,}{f(Y)}\mathbf{1}_{\{f(Y)\leq e_{\mathbf{n}}\}}\right\Vert _{2}=\mathcal{O}\left(\frac{1}{\hat{\mathbf{n}}^{\frac{1}{2}+\delta}}\right)$
by $\mathbf{E}\left(\exp\left(\left\Vert r(Y)\right\Vert \,\mathbf{1}_{\{f(Y)\leq e_{\mathbf{n}}\}}\right)\right)=\mathcal{O}\left(\hat{\mathbf{n}}^{-\xi}\right)$
for some $\xi>0$. Then, \begin{eqnarray*}
P\left(\left\Vert \frac{1}{\hat{\mathbf{n}}}\sum_{\mathbf{i}\in\mathcal{I}_{\mathbf{n}}}\frac{r(Y_{\mathbf{i}})\,}{f(Y_{\mathbf{i}})}\mathbf{1}_{\{f(Y_{\mathbf{i}})\leq e_{\mathbf{n}}\}}\right\Vert >\frac{\varepsilon}{v_{\mathbf{n}}}\right) & \leq & P\left(\frac{C}{\hat{\mathbf{n}}}\sum_{\mathbf{i}\in\mathcal{I}_{\mathbf{n}}}\left\Vert r(Y_{\mathbf{i}})\right\Vert \,\mathbf{1}_{\{f(Y_{\mathbf{i}})\leq e_{\mathbf{n}}\}}>\frac{\varepsilon}{v_{\mathbf{n}}}\right)\end{eqnarray*}
and because of Minskovski's inequality: for all $k\in\mathbb{N}^{*}$,
$\mathbf{E}\left(\left(\frac{1}{\hat{\mathbf{n}}}\sum_{\mathbf{i}\in\mathcal{I}_{\mathbf{n}}}\left\Vert r(Y_{\mathbf{i}})\right\Vert \,\mathbf{1}_{\{f(Y_{\mathbf{i}})\leq e_{\mathbf{n}}\}}\right)^{k}\right)\leq\left\Vert r(Y)\,\mathbf{1}_{\{f(Y)\leq e_{\mathbf{n}}\}}\right\Vert _{k}^{k}\,$,
we can say that with using the argument $\mathbf{E}\left(\exp\left(\left\Vert r(Y)\right\Vert \,\mathbf{1}_{\{f(Y)\leq e_{\mathbf{n}}\}}\right)\right)=\mathcal{O}\left(\hat{\mathbf{n}}^{-\xi}\right)$
: \begin{eqnarray*}
P\left(\left\Vert \frac{1}{\hat{\mathbf{n}}}\sum_{\mathbf{i}\in\mathcal{I}_{\mathbf{n}}}\frac{r(Y_{\mathbf{i}})\,}{f(Y_{\mathbf{i}})}\mathbf{1}_{\{f(Y_{\mathbf{i}})\leq e_{\mathbf{n}}\}}\right\Vert >\frac{\varepsilon}{v_{\mathbf{n}}}\right) & \leq & \mathbf{E}\,\left[\exp\left(\left\Vert r(Y)\right\Vert \,\mathbf{1}_{\{f(Y)\leq e_{\mathbf{n}}\}}\right)\right]\exp\left(-\frac{\varepsilon}{v_{\mathbf{n}}}\right)\\
 & \leq & C_{1}\,\hat{\mathbf{n}}^{-\xi}.\exp\left(-\varepsilon\left(\frac{\log\log\,\hat{\mathbf{n}}}{\hat{\mathbf{n}}}\right)^{-\frac{1}{2}}\right)\,\,\mbox{for some}\, C_{1}>0.\\
 & \leq & C_{1}\,\exp\left(-\xi\,\log\hat{\mathbf{n}}\,-\varepsilon\left(\frac{\hat{\mathbf{n}}}{\log\log\,\hat{\mathbf{n}}}\right)^{\frac{1}{2}}\right)\\
 & \leq & C_{1}\,\exp\left(-\min(\xi\,,\varepsilon)\,\left(\log\hat{\mathbf{n}}\,+\frac{\sqrt{\log\hat{\mathbf{n}}}}{\sqrt{\log\log\,\hat{\mathbf{n}}}}\right)\right)\\
 & \leq & C_{1}\,\exp\left(-\min(\xi\,,\varepsilon)\,\log\hat{\mathbf{n}}\,\left(1+\frac{1}{\sqrt{\left(\log\hat{\mathbf{n}}\right)\log\log\,\hat{\mathbf{n}}}}\right)\right)\end{eqnarray*}
as $\hat{\mathbf{n}}\to+\infty$, $\exp\left(-\min(\xi\,,\varepsilon)\,\log\hat{\mathbf{n}}\,\left(1+\frac{1}{\sqrt{\left(\log\hat{\mathbf{n}}\right)\log\log\,\hat{\mathbf{n}}}}\right)\right)\simeq\hat{\mathbf{n}}^{-C_{2}}$
where $c_{2}$ is positive constant. So, $\frac{1}{\hat{\mathbf{n}}}\sum_{\mathbf{i}\in\mathcal{I}_{\mathbf{n}}}\frac{r(Y_{\mathbf{i}})\,}{f(Y_{\mathbf{i}})}\mathbf{1}_{\{f(Y_{\mathbf{i}})\leq e_{\mathbf{n}}\}}=o_{a.s}\left(\left(\frac{\log\log\,\hat{\mathbf{n}}}{\hat{\mathbf{n}}}\right)^{\frac{1}{2}}\right)$
and the proof is complet by using Lemma \ref{thm:Carbon-et-alAS}
and sketching the proof of Theorem \ref{thm:conv-prob}~~~~$\square$

\subsection*{Proof of Corollary \ref{cor:Conv-a.s-opt-rate} }

If moreover we chose $h_{\mathbf{n}}\simeq\hat{\mathbf{n}}^{-c_{1}}$
and $e_{\mathbf{n}}\simeq\hat{\mathbf{n}}^{-c_{2}}$ where $\frac{c_{2}}{k}+\frac{1}{4k}<c_{1}<\frac{1}{2}-2c_{2}$,
then, \[
\sqrt{\frac{\hat{\mathbf{n}}}{\log\log\hat{\mathbf{n}}}}\times\frac{\Psi_{\mathbf{n}}^{2}}{e_{\mathbf{n}}^{2}}=\sqrt{\frac{\hat{\mathbf{n}}}{\log\log\hat{\mathbf{n}}}}\times\left(\hat{\mathbf{n}}^{-2kc_{1}+2c_{2}}+\hat{\mathbf{n}}^{-1+c_{1}+2c_{2}}\log\hat{\mathbf{n}}\right)\]
$\sqrt{\frac{\hat{\mathbf{n}}}{\log\log\hat{\mathbf{n}}}}\times\frac{\Psi_{\mathbf{n}}^{2}}{e_{\mathbf{n}}^{2}}=\frac{\hat{\mathbf{n}}^{\frac{1}{2}-2kc_{1}+2c_{2}}}{\sqrt{\log\log\hat{\mathbf{n}}}}+\hat{\mathbf{n}}^{-\frac{1}{2}+c_{1}+2c_{2}}\log\hat{\mathbf{n}}$
this latter tend to zero as soon as $\frac{c_{2}}{k}+\frac{1}{4k}\leq c_{1}<\frac{1}{2}-2c_{2}$
\\
 The proof is obtained by sketching the proof of Corolary \ref{cor:cov-prob-vitesse opt}
and using the law of the iterated logarithm recalled in Lemma \ref{lem: gLaw-of-large}.

\bibliographystyle{plain}

\bibliographystyle{ref1} \bibliographystyle{ref1}
\bibliography{ref1}

\end{document}